%% file: main.tex
\documentclass[onecolumn,10pt,oneside,reqno]{article}

\input{header}

\newcommand{\RefereeOne}[1]{\textcolor{black}{#1}}
\newcommand{\RefereeTwo}[1]{\textcolor{black}{#1}}
\newcommand{\RefereeALL}[1]{\textcolor{black}{#1}}

\begin{document}

\maketitle

%
%
%

\begin{abstract}
Data-driven transformations that reformulate nonlinear systems in a linear framework have the potential to enable the prediction, estimation, and control of strongly nonlinear dynamics using linear systems theory.  
The Koopman operator has emerged as a principled linear embedding of nonlinear dynamics, and its eigenfunctions establish intrinsic coordinates along which the dynamics behave linearly.  
Previous studies have used finite-dimensional approximations of the Koopman operator for model-predictive control approaches. In this work, we illustrate a fundamental closure issue of this approach and argue that it is beneficial to
\RefereeOne{first validate eigenfunctions and then construct reduced-order models in these validated eigenfunctions}. 
These coordinates form a Koopman-invariant subspace \textit{by design} and, thus, have improved predictive power.
\RefereeOne{We show then how the control can be formulated directly in these intrinsic coordinates and discuss potential benefits and caveats of this perspective. The resulting control architecture is termed \textit{Koopman Reduced Order Nonlinear Identification and Control} (KRONIC).}
It is further demonstrated that these eigenfunctions can be approximated with data-driven regression and power series expansions, based on the partial differential equation governing the infinitesimal generator of the Koopman operator.
Validating discovered eigenfunctions is crucial and we show that lightly damped eigenfunctions may be faithfully extracted 
\RefereeOne{from EDMD or an implicit formulation}.
These lightly damped eigenfunctions are particularly relevant for control, as they correspond to nearly conserved quantities that are associated with persistent dynamics, such as the Hamiltonian.  
KRONIC is then demonstrated on a number of relevant examples, including 1) a nonlinear system with a known linear embedding, 2) a variety of Hamiltonian systems, and 3) a high-dimensional double-gyre model for ocean mixing.  \\
\end{abstract}
\hspace{1cm}{\it Keywords\/}: Dynamical systems, nonlinear control, Koopman theory, system identification, machine learning.


\section{Introduction}\label{Sec:Introduction}
\input{Introduction}

\section{Motivation}\label{Sec:Motivation}  
\input{Motivation}

\section{Background}\label{Sec:Background}  
\input{Background}

\section{Koopman operator control in eigenfunctions}\label{Sec:KRONIC}
\input{KRONIC}

\section{Identifying Koopman eigenfunctions from data}\label{Sec:Regression}
\input{Regression}

\section{Example: System with a slow manifold}\label{Sec:SlowManifold}
\input{SlowManifold}

\section{Example: Hamiltonian energy control}\label{Sec:HamiltonianSystems}
\input{HamiltonianSystems}

\section{Example: Basin hopping in a double well}\label{Sec:BasinHopping}
\input{AsymmetricPotentialWell}

\section{Example: Double Gyre flow}\label{Sec:DoubleGyre}
\input{DoubleGyre}

\section{Discussion and conclusions}\label{Sec:Discussion}
\input{Conclusions}

\section*{Acknowledgments}
EK acknowledges funding by the Moore/Sloan foundation, the Washington Research Foundation, and the eScience Institute.  
SLB and JNK acknowledge support from the Defense Advanced Research Projects Agency (DARPA contract HR011-16-C-0016) and the UW Engineering Data Science Institute, NSF HDR award \#1934292.  
SLB acknowledges support from the Army Research Office (W911NF-17-1-0306) and the Air Force Office of Scientific Research (FA9550-16-1-0650). 
JNK acknowledges support from the Air Force Office of Scientific Research (FA9550-15-1-0385). 
The authors gratefully acknowledge many valuable discussions with Josh Proctor about Koopman theory and extensions to control.  
We would also like to acknowledge Igor Mezi\'c, Maria Fonoberova, Bernd Noack, Clancy Rowley, Sam Taira, and Lionel Mathelin.

\setcounter{section}{0}
\renewcommand\thesection{\Alph{section}}
\section{Optimal control}\label{Sec:Control}
\input{ControlAppendix}
\section{Effect of misrepresentation of eigenfunctions}\label{Sec:ControlError}
\input{ControlError}
\section{Series solutions for eigenfunctions}\label{Sec:AnalyticalKoopmanEfun}
\input{AnalyticalEigenfunction}

\bibliographystyle{plain} 
\bibliography{references}

\end{document}

%% file: header.tex
\usepackage[top=1in, bottom=1in, left=.85in, right=.85in,paperwidth=8.5in,paperheight=11in]{geometry}

\usepackage[svgnames]{xcolor}
\usepackage{amsmath}
\usepackage{amssymb}
\usepackage{amsfonts}
\usepackage{amsopn}
\usepackage{graphicx}
\usepackage{epstopdf}
\usepackage{overpic}
\usepackage{url}
\usepackage{caption}
\usepackage{rotating}
\usepackage{multirow}
\usepackage{mathtools}
\usepackage[innercaption]{sidecap}
\usepackage{setspace}
\usepackage{hyperref}
\usepackage{psfrag}
\usepackage{lipsum}

\usepackage{algorithmic} 
\usepackage{algorithm}

\usepackage{enumitem}
\setlist[enumerate]{leftmargin=.5in}
\setlist[itemize]{leftmargin=.5in}

\newcommand{\bb}{\mathbf{b}}
\newcommand{\bc}{\mathbf{c}}

\newcommand{\bq}{\mathbf{q}}
\newcommand{\bp}{\mathbf{p}}

\newcommand{\bu}{\mathbf{u}}

\newcommand{\bx}{\mathbf{x}}

\newcommand{\by}{\mathbf{y}}
\newcommand{\bz}{\mathbf{z}}
\newcommand{\bA}{\mathbf{A}}
\newcommand{\bB}{\mathbf{B}}
\newcommand{\bC}{\mathbf{C}}
\newcommand{\bD}{\mathbf{D}}
\newcommand{\bK}{\mathbf{K}}

\newcommand{\bP}{\mathbf{P}}
\newcommand{\bQ}{\mathbf{Q}}
\newcommand{\bR}{\mathbf{R}}

\newcommand{\bU}{\mathbf{U}}
\newcommand{\bV}{\mathbf{V}}

\newcommand{\bX}{\mathbf{X}}

\newcommand{\bxi}{\boldsymbol{\xi}}

\newcommand{\bvarphi}{\boldsymbol{\varphi}}

\newcommand{\bSigma}{\boldsymbol{\Sigma}}
\newcommand{\bTheta}{\boldsymbol{\Theta}}
\newcommand{\bLambda}{\boldsymbol{\Lambda}}
\newcommand{\bGamma}{\boldsymbol{\Gamma}}

\ifpdf
\DeclareGraphicsExtensions{.eps,.pdf,.png,.jpg}
\else
\DeclareGraphicsExtensions{.eps}
\fi

\title{\huge{\textbf{Data-driven discovery of Koopman eigenfunctions for control}}}
\author{Eurika Kaiser$^{1}$\footnote{E-mail: \href{mailto:eurika@uw.edu}{eurika@uw.edu}},
	J. Nathan Kutz$^{2}$, and Steven L. Brunton$^{1}$\\~\\
\small{$^1$ Department of Mechanical Engineering, University of Washington, Seattle, WA 98195, United States}\\
\small{$^2$ Department of Applied Mathematics, University of Washington, Seattle, WA 98195,	United States}
\\~\\
}
\date{}
\date{\today}

%% file: Introduction.tex
In contrast to linear systems, a generally applicable and scalable framework for the control of nonlinear systems remains an engineering grand challenge.
Improved nonlinear control has the potential to transform our ability to interact with and manipulate complex systems across broad scientific, technological, and industrial domains. 
From turbulence control to brain-machine interfaces, emerging technologies are characterized by high-dimensional, strongly nonlinear, and multiscale phenomena that lack simple models suitable for control design.  
This lack of simple equations motivates \emph{data-driven} control methodologies, which include system identification for model discovery~\cite{ljung:book,Bamieh2002ijrnc,Nelles2013book,Billings2013book,Brunton2016pnas}.  
Alternatively, one can seek transformations that embed nonlinear dynamics in a global linear representation, as in the Koopman framework~\cite{Koopman1931pnas,Mezic2005nd}.
The goal of this work is to reformulate controlled nonlinear dynamics in a Koopman-eigenfunction framework, 
referred to as \emph{Koopman Reduced Order Nonlinear Identification and Control} (KRONIC),  
that shows improved predictive power and is amenable to powerful linear optimal and robust control techniques~\cite{sp:book,dp:book,stengel2012book}.

A wide range of data-driven and nonlinear control approaches exist in the literature, including model-free adaptive control~\cite{Krstic:1995}, extremum-seeking~\cite{Krstic2003book}, 
gain scheduling~\cite{Rugh2000automatica}, feedback linearization~\cite{charlet1989dynamic}, describing functions~\cite{vander1968multiple}, sliding mode control~\cite{edwards1998sliding}, singular perturbation~\cite{kokotovic1976singular}, geometric control~\cite{Brockett1976automatica}, back-stepping~\cite{kokotovic1992joy}, model predictive control~\cite{camacho2013model,mayne2000constrained}, 
reinforcement learning~\cite{Sutton1998book}, 
and machine learning control~\cite{Hansen2009ieeetec,Brunton2015amr}. 
Although considerable progress has been made in the control of nonlinear systems~\cite{Krstic:1995,isidori2013nonlinear,sastry2013nonlinear}, methods are generally tailored to a specific class of problems, require considerable mathematical and computational resources, or don't readily generalize to new applications.  
Currently there is no overarching framework for nonlinear control, that is generically applicable to a wide class of potentially high-dimensional systems, as exists for linear systems~\cite{dp:book,sp:book}.  
Generally applicable frameworks, such as dynamic programming~\cite{Bertsekas2005book}, Pontryagin's maximum principle~\cite{Pontryagin2062interscience}, and model-predictive control (MPC)~\cite{allgower1999springer,camacho2013model}, often suffer from the curse of dimensionality and require considerably computational effort, e.g. solving adjoint equations, when applied to nonlinear systems, which can be mitigated to some degree by combining these with low-dimensional, linear representations. 
Fortunately, the rise of big data, advances in machine learning, and new approaches in dynamical systems are changing how we approach these canonically challenging nonlinear control problems.  
For instance, recently deep reinforcement learning has been combined with MPC~\cite{Peng.2009,Zhang2016icra}, yielding impressive results in the large-data limit.

\begin{figure}
	\centering
	\includegraphics[width=0.6\textwidth]{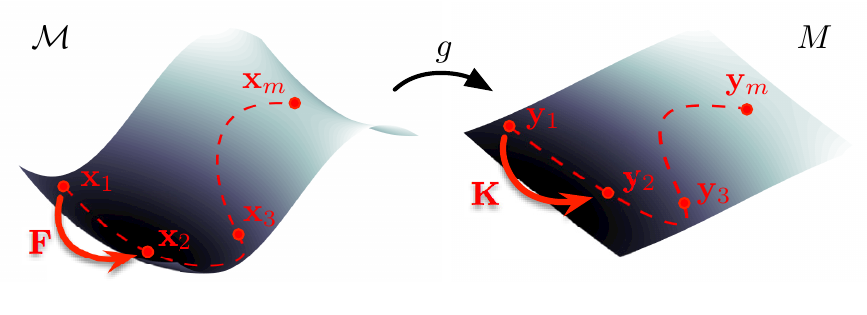}
	\vspace{-.2in}	
	\caption{Koopman embedding to linearize nonlinear dynamics.}\label{Fig:KoopmanEmbedding}
	\vspace{-.25in}	
	\label{Fig:Schematic:Embedding}
\end{figure}

Koopman operator theory has recently emerged as a leading framework to obtain linear representations of nonlinear dynamical systems from data~\cite{Mezic2005nd}.  
This operator-theoretic perspective complements the more standard geometric~\cite{guckenheimer_holmes} and probabilistic~\cite{Dellnitz2001book,Dellnitz2002hds} perspectives.  
The ability to embed nonlinear dynamics in a linear framework 
(see Fig.~\ref{Fig:Schematic:Embedding}) 
is particularly promising for the prediction, estimation, and control of nonlinear systems.  
In 1931, Koopman showed that a nonlinear dynamical system may be represented by an infinite-dimensional linear operator acting on the space of measurement functions of the state of the system~\cite{Koopman1931pnas}.  
Formulating dynamics in terms of measurements is appealing in the era of big data.  
Since the seminal work of Mezi\'c and Banaszuk~\cite{Mezic2004physicad} and Mezi\'c~\cite{Mezic2005nd}, Koopman theory has been the focus of efforts to characterize nonlinear systems.    
Many classical results have been extended to the Koopman formalism~\cite{Budivsic2012chaos,Mezic2013arfm}.  
For example, level sets of Koopman eigenfunctions form invariant partitions~\cite{Budivsic2012physd} and may be used to analyze mixing.  
The Hartman-Grobman theorem has also been generalized to provide a linearizing  transform in the entire basin of attraction of a stable or unstable equilibrium or periodic orbit~\cite{Lan2013physd}.
 
Recently, Koopman theory has been applied for system identification~\cite{Mezic2015cdc,Proctor2016arxiv,williams2016ifac,mauroy2016cdc},  estimation~\cite{Surana2016cdc,Surana2016nolcos} and control~\cite{wilson2014jads,Brunton2016plosone,Korda2016arxiv,sootla2017arxiv,abraham2017arxiv} of nonlinear systems.
The Koopman operator is infinite-dimensional, and control laws are typically based on a finite-dimensional approximation.  
Dynamic mode decomposition (DMD)~\cite{Schmid2010jfm,Rowley2009jfm,Tu2014jcd,Kutz2016book} approximates the Koopman operator with a best-fit linear model.  
However, DMD is based on linear measurements, which do not span a Koopman invariant subspace for many nonlinear systems~\cite{Williams2015jnls,Brunton2016plosone,Kutz2016book}.  
Current data-driven methods to approximate the Koopman operator include extended DMD (EDMD)~\cite{Williams2015jnls,Williams2015jcd,klus2015numerical} 
and the variational approach of conformation dynamics (VAC)~\cite{noe2013variational,nuske2016variational}. 
EDMD was recently used for model predictive control with promising results~\cite{Korda2016arxiv}.  
However, EDMD models may suffer from closure issues for systems with multiple fixed points or attractors, as a linear model only has a single fixed point~\cite{Brunton2016plosone},
which may lead to corrupted dynamics~\cite{Brunton2016plosone,Kutz2016book} and emphasizes the importance of model validation. 
For instance, some eigenfunctions may be distorted when projected onto a finite-dimensional measurement subspace, and it may be advantageous to construct a reduced-order description from {\it validated} eigenfunctions, that exhibit behavior as predicted by their associated eigenvalue. 
For chaotic systems, delay coordinates provides a promising embedding~\cite{Susuki2015cdc,Brunton2017natcomm,Arbabi2016arxiv}. 
Obtaining useful data-driven coordinate transformations that approximate Koopman eigenfunctions remains an open challenge in data-driven dynamical systems~\cite{Kutz2016book,Brunton2016plosone}.  

In the present work, we 
\RefereeOne{build on the existing EDMD and}
reformulate the Koopman-based control problem in eigenfunction coordinates and provide strategies to identify lightly damped eigenfunctions from data, that can be subsequently used for control. In particular:
\begin{itemize}
	\item Koopman eigenfunctions provide a principled linear embedding of nonlinear dynamics resulting in an intrinsic coordinate system, in which the system is closed under the action of the Koopman operator. 
	We formulate the incorporation of control in these intrinsic coordinates and discuss benefits and caveats of this approach.
	Further, Koopman eigenfunctions can be associated with geometric system properties and coherent structures. 
	Thus, the designed controller may be employed to manipulate particular coherent structures. 
	For example, the Hamiltonian energy is a Koopman eigenfunction,  and we are able to control the system by manipulating this function.  
	\item Smooth eigenfunctions in the point spectrum of the Koopman operator can be discovered from given data using sparse regression providing interpretable representations. 
	We propose sparsity-promoting algorithms to regularize EDMD or to discover eigenfunctions directly in an implicit formulation.
	\item We further demonstrate the importance of model validation to distinguish {\it accurately identified} eigenfunctions from {\it spurious} ones. Lightly damped eigenfunctions are often not corrupted and can be used to construct reduced-order Koopman models.
\end{itemize}

\begin{figure*}
	\vspace{.1in}
	\begin{overpic}[width=1\textwidth]{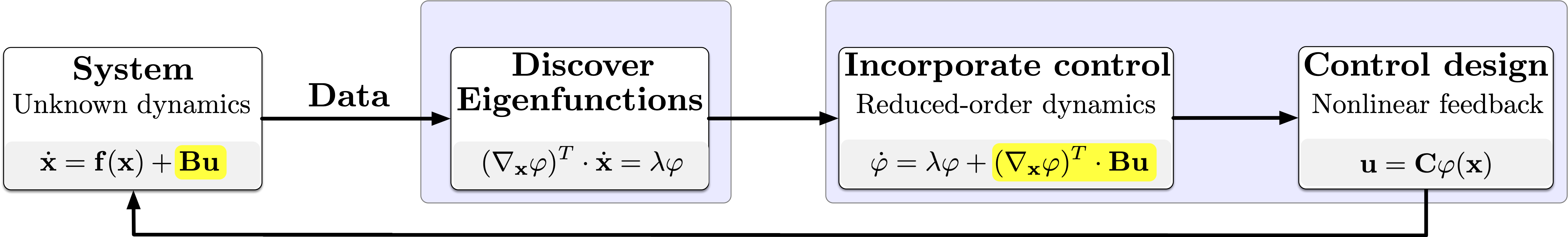} 
		\put(32.5,12.5){Section~\ref{Sec:Regression}}
		\put(73.5,12.5){Section~\ref{Sec:KRONIC}}
	\end{overpic}	
	\caption{Control of nonlinear systems via reduced Koopman-invariant representations in eigenfunction coordinates.}
	\label{Fig:Schematic:KRONIC}
	\vspace{-.05in}
\end{figure*}

These nonlinear control techniques generalize to \emph{any} lightly damped eigenfunction.  
As a more sophisticated example, we consider the double gyre flow, which is a model for ocean mixing.  
The discovery of intrinsic coordinates for optimized nonlinear control establishes our data-driven 
KRONIC framework\footnote{Code at \url{https://github.com/eurika-kaiser/KRONIC}.}, shown in Fig.~\ref{Fig:Schematic:KRONIC}.  

The present work is outlined as follows: In Sec.~\ref{Sec:Motivation}, we demonstrate the importance of eigenfunction validation and motivate the use of sparse regression for their discovery. In Sec.~\ref{Sec:Background}, key results in Koopman spectral theory and corresponding data-driven approaches are summarized, and a brief background on optimal control is provided. 
The approach for identifying of Koopman eigenfunctions from data using sparse regression is outlined in Sec.~\ref{Sec:Regression}. In Sec.~\ref{Sec:KRONIC}, it is shown how control can be incorporated in the eigenfunction formulation. 
An analytical example is examined in Sec.~\ref{Sec:SlowManifold} to illustrate the control problem in terms of Koopman eigenfunction coordinates. 
The KRONIC framework is then demonstrated on several Hamiltonian systems (Sec.~\ref{Sec:HamiltonianSystems}), for basin-hopping in an asymmetric double potential well (Sec.~\ref{Sec:BasinHopping}), and the autonomous and non-autonomous double gyre flow (Sec.~\ref{Sec:DoubleGyre}).
A discussion and outlook on future directions is provided in Sec.~\ref{Sec:Discussion}.

%% file: Motivation.tex
\RefereeOne{
Finite-dimensional approximations of the Koopman operator are typically obtained as its projection onto a specified basis or dictionary.
Extended dynamic mode decomposition (EDMD)~\cite{williams2016ifac} has emerged as the leading numerical approach by solving a least-squares problem.
A well-known issue arises when trying to identify the full operator in a finite set of basis functions, that sometimes the finite-dimensional approximation is not closed and spurious eigenfunctions may appear.
For this reason, it can be important to perform consistency checks, such as validating the linearity property of eigenfunctions. 
The present work builds on EDMD addressing this limitation by re-formulating the regression problem for the direct identification of Koopman eigenfunctions.
Further, we demonstrate how EDMD may be regularized to obtain more accurate eigenfunction representations. 
In the following, we illustrate the closure problem using a polynomial basis.
}
As a motivating example (examined in detail in Sec.~\ref{Sec:SlowManifold}), we consider a system with quadratic nonlinearity that gives rise to a slow manifold~\cite{Tu2014jcd}:
\begin{align}\label{Eqn:Motivation:SlowManifold}
	\frac{\mathrm d}{\mathrm dt}\begin{bmatrix} x_1 \\ x_2 \end{bmatrix} = \begin{bmatrix} \mu x_1 \\ \lambda (x_2 - x_1^2)\end{bmatrix}
\end{align}
with $\mu=-0.1$ and $\lambda=-1$.
By a clever choice of observable functions $\by$, the nonlinear system~\eqref{Eqn:Motivation:SlowManifold} may be represented as a linear system: 
\begin{equation}
\frac{d}{dt}{\bf y} = {\bf K} {\bf y},
\end{equation} 
where ${\bf K}$ represents a finite-dimensional approximation of the Koopman operator.
The system~\eqref{Eqn:Motivation:SlowManifold} is one of few analytical examples, for which a closed, finite-dimensional, linear Koopman approximation exists.

\RefereeALL{
The EDMD model fit on the first $9$ monomials (up to third degree) is shown in Fig.~\ref{Fig:EDMD_EvecPrediction1}(a). 
\begin{figure}[tb]
	\centering
	\includegraphics[width=\textwidth]{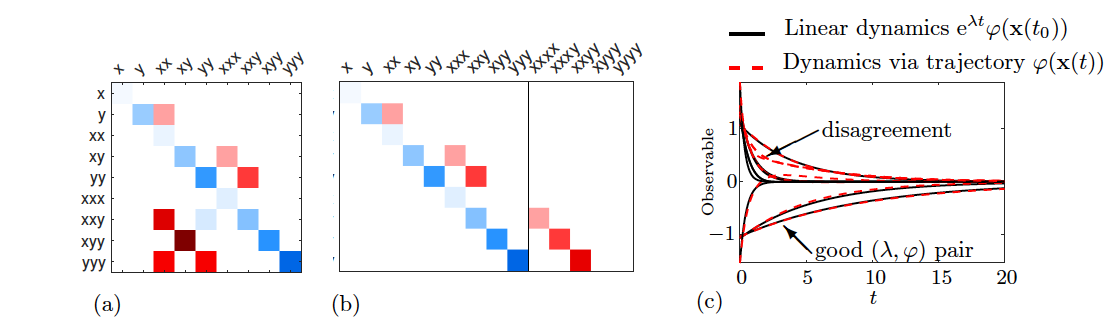}
	\caption{
	\RefereeOne{	
		EDMD models give rise to inaccurate eigenvectors, demonstrated for the slow manifold system: 
		(a) EDMD regression model on up to third degree monomials, 
		(b) fitting onto higher-degree monomials, 
		(c) evolution of eigenvectors identified from the model in (a), which are evaluated on a single test trajectory, and as predicted using the corresponding eigenvalue.
	}
	}
	\vspace{-.08in}
	\label{Fig:EDMD_EvecPrediction1}
\end{figure}
Some of the eigenvectors are \emph{spurious} (see Fig.~\ref{Fig:EDMD_EvecPrediction1}(c)), i.e. 
\RefereeOne{
the evolution of the eigenfunction $\varphi(\bx(t))$ obtained by evaluating the eigenfunction on a trajectory $\bx(t)$ does not correspond to the linear prediction using the eigenvalue, $\mathrm{e}^{\lambda t}\varphi(\bx(t_0))$.}	
The reason is that the model is not closed in the particular choice of basis functions, because the derivatives of third degree monomials give rise to higher-degree monomials (see Fig.~\ref{Fig:EDMD_EvecPrediction1}(b)) resulting in a nonlinear model:
\begin{align}
	\frac{\mathrm d}{\mathrm dt}\begin{bmatrix} y_1 \\ y_2 \\ y_3\\ y_4\\ y_5\\y_6\\y_7\\y_8\\y_9\end{bmatrix} = \begin{bmatrix} 
		\mu & 0 & 0 & 0 & 0 & 0 & 0 & 0 & 0\\ 
		0 & \lambda & -\lambda & 0 & 0 & 0 & 0 & 0 & 0\\ 
		0 & 0 &  2\mu & 0 & 0& 0 & 0 & 0 & 0\\ 
		0 & 0 & 0 & \mu+\lambda & -\lambda & 0 & 0 & 0 & 0\\ 
		0 & 0 & 0 & 0 & 2\lambda & 0 & -2\lambda & 0 & 0\\
		0 & 0 & 0 & 0 & 0 & 3\mu & 0 & 0 & 0 \\
		0 & 0 & 0 & 0 & 0 & 0       & \lambda+2\mu & 0 & 0\\
		0 & 0 & 0 & 0 & 0 & 0 & 0 & 2\lambda+\mu & 0\\
		0 & 0 & 0 & 0 & 0 & 0 & 0 & 0 & 3\lambda\\
	\end{bmatrix}
	\begin{bmatrix}
		y_1 \\ y_2\\ y_3\\y_4\\y_5\\y_6\\y_7\\y_8\\y_9
	\end{bmatrix}
	+ 
	\begin{bmatrix}
		0\\0\\0\\0\\0\\0\\-\lambda x_1^4\\-2\lambda x_1^3x_2\\-3\lambda x_1^2x_2^2
	\end{bmatrix},
\end{align}
where $(y_1,y_2,y_3,y_4,y_5, y_6, y_7, y_8,y_9) = (x_1,x_2,x_1^2,x_1x_2,x_2^2, x_1^3, x_1^2x_2,x_1x_2^2,x_3^3)$.
These additional terms, that are not in the span of $\{y_1,\cdots,y_9\}$, will be \emph{aliased} in the corresponding row equations corrupting the system matrix $\bK$.  Thus, some of the eigenfunctions will be spurious, affecting the prediction accuracy of the model based on $\bK$.  
However, it may be possible to identify a subset of the eigenfunctions that are not corrupted, e.g. those eigenpairs that show good agreement in Fig.~\ref{Fig:EDMD_EvecPrediction1}(c), and use these to construct a reduced-order model with improved prediction.
Alternatively, EDMD may be regularized using sparsity-promoting techniques to regress a (approximate) closed model on a subset of basis functions.
For instance, choosing the five observable functions $(y_1,y_2,y_3,y_4,y_5) = (x_1,x_2,x_1^2,x_1x_2,x_1^3)$, which are a subset of the nine monomials used above, yields:
\begin{align}\label{Eqn:Motivation:SlowManifoldKoopman5D}
\frac{\mathrm d}{\mathrm dt}\begin{bmatrix} y_1 \\ y_2 \\ y_3\\ y_4\\ y_5\end{bmatrix} = \begin{bmatrix} \mu & 0 & 0 & 0 & 0\\ 0 & \lambda & -\lambda & 0 & 0\\ 0 & 0 &  2\mu & 0 & 0\\ 0 & 0 & 0 & \mu+\lambda & -\lambda\\ 0 & 0 & 0 & 0 & 3\mu\end{bmatrix}
\begin{bmatrix}
y_1 \\ y_2\\ y_3\\y_4\\y_5
\end{bmatrix},
\end{align}
which is a 5-dim. linear system, that remains closed under the action of the Koopman operator.
}

\RefereeALL{
Measurements from real-world systems are generally corrupted by noise, which can be more challenging for model identification procedures.
Figure~\ref{Fig:EDMD_EvecPrediction1_Noise}(c) shows the poor prediction performance of the EDMD model trained on noisy data (displayed in Fig.~\ref{Fig:EDMD_EvecPrediction1_Noise}(a)).
\begin{figure}[tb]
\centering
\includegraphics[width=\textwidth]{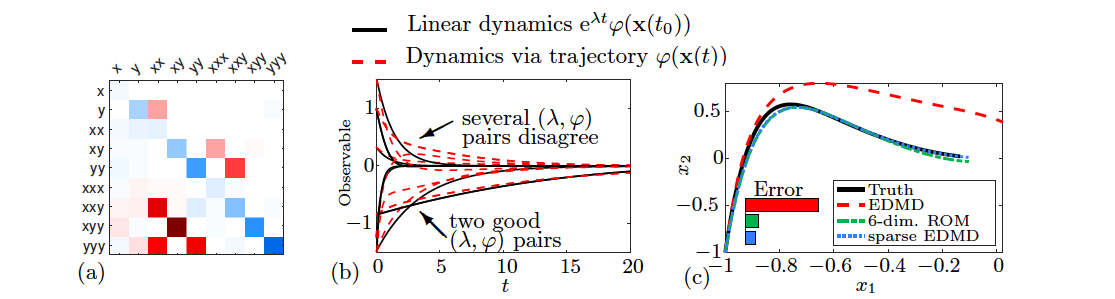}
\caption{
\RefereeOne{
Sparsified EDMD models increase robustness and prevent overfitting on noisy data:
(a) EDMD model,
(b) validation of eigenpairs, and
(c) phase plot with $L_2$ error for different models. 
For details see Fig.~\ref{Fig:EDMD_EvecPrediction1} and text.}
}
\label{Fig:EDMD_EvecPrediction1_Noise}
\end{figure}
More eigenfunctions are inaccurate as a result of an increasing number of non-vanishing coefficients in the state matrix compared with the noise-free situation in Fig.~\ref{Fig:EDMD_EvecPrediction1}.
The least-squares solution overfits resulting in a full matrix with small, but non-vanishing coefficients.
By sparsifying the EDMD state matrix or constructing a reduced-order model based on accurate eigenfunctions (validated from Fig.~\ref{Fig:EDMD_EvecPrediction1_Noise}(b)) higher prediction accuracy and robustness to noise can be achieved (compare models in Figure~\ref{Fig:EDMD_EvecPrediction1_Noise}(c)).
Ideally, it would be possible to learn \emph{good} eigenfunctions directly, which are designed to behave linearly and evolve as predicted by their associated eigenvalue, which would potentially significantly increase prediction accuracy and reduce the dimension of the model.
}

%% file: Background.tex
In this section, we provide a brief background on Koopman spectral theory in Sec.~\ref{Sec:Background:Koopman}
and a numerical algorithm, dynamic mode decomposition, to approximate the Koopman operator in Sec.~\ref{Sec:Background:DMD}. 
Key results in optimal control are then summarized in Sec.~\ref{Sec:Background:OptimalControl}.

\subsection{Koopman spectral theory}\label{Sec:Background:Koopman}

The classical geometric theory of dynamical systems considers a set of coupled ordinary differential equations
\begin{equation}\label{Eq:Dynamics} 
\frac{\mathrm d}{\mathrm dt}\bx(t) =  {\boldsymbol{\it f}}(\bx)
\end{equation}
in terms of the state of the system $\bx\in\mathcal{M}$, where $\mathcal{M}$ is a differentiable manifold, often given by $\mathcal{M}=\mathbb{R}^n$.  
In discrete time, the dynamics are given by 
\begin{equation}
\bx_{k+1} = \mathbf{F}(\bx_k),\label{Eq:DiscreteDynamics}
\end{equation}
where $\mathbf{F}$ may be the flow map of the dynamics in \eqref{Eq:Dynamics}:
\begin{equation}
\mathbf{F}(\bx(t_0)) = \bx(t_0)+\int_{t_0}^{t_0+t} {\boldsymbol{\it f}}(\bx(\tau))\,d\tau.
\end{equation}
Discrete-time systems are more general and form a superset, containing those induced by continuous-time dynamics.  
Moreover, discrete-time dynamics are often more consistent with experimental measurements, and may be preferred for numerical analysis.  
The geometric perspective then considers fixed points and invariant structures of the dynamics.  

In 1931, B. O. Koopman introduced the operator theoretic perspective, showing that there exists an infinite-dimensional linear operator, given by $\mathcal{K}$, that acts to advance all measurement functions $g: \mathcal{M}\rightarrow \mathbb{R}$ of the state with the flow of the dynamics:  
\begin{equation}
Kg = g\circ\mathbf{F}.\label{Eq:KoopmanDiscrete}
\end{equation}
Thus, the Koopman operator advances measurements linearly: 
\begin{equation}
g(\bx_{k+1}) = K g(\bx_k).\label{Eq:KoopmanDiscrete2}
\end{equation}
For smooth dynamics, there is a continuous system
\RefereeTwo{
\begin{equation}
\frac{\mathrm d}{\mathrm dt} g(\bx) = \mathcal{K}g(\bx) = \nabla g(\bx)\cdot {\bf f}(\bx),
\end{equation}
}
where $\mathcal{K}$ is the infinitesimal generator of the one-parameter family of Koopman operators $K$.  

The Koopman operator is linear, which is appealing, but is infinite dimensional, posing issues for representation and computation.  
Instead of capturing the evolution of all measurement functions in a Hilbert space, applied Koopman analysis approximates the evolution on a subspace spanned by a finite set of measurement functions.  
It is possible to obtain a finite-dimensional matrix representation of the Koopman operator by restricting it to an invariant subspace.  
A Koopman invariant subspace is spanned by any set of eigenfunctions of the Koopman operator.  
A Koopman eigenfunction $\varphi(\bx)$ corresponding to eigenvalue $\lambda$ satisfies
\begin{equation}
\lambda\varphi(\bx) = \varphi(\mathbf{F}(\bx)).
\end{equation}
In continuous-time, a Koopman eigenfunction $\varphi(\bx): \mathcal{M}\rightarrow \mathbb{C}$ satisfies
\begin{equation}
\frac{\mathrm d}{\mathrm dt}\varphi(\bx) = \lambda \varphi(\bx).\label{Eq:KoopmanEfun}
\end{equation}
Obtaining Koopman eigenfunctions from data or analytically is a central applied challenge in modern dynamical systems.  
Discovering these eigenfunctions enables globally linear representations of strongly nonlinear systems in terms of these \emph{intrinsic} observables.  
\RefereeTwo{The evolution equation~\eqref{Eq:KoopmanEfun} describes the unactuated behavior, which will be extended to incorporate the effect of control in our KRONIC framework (compare the third box in Fig.~\ref{Fig:Schematic:KRONIC} and for further details we refer to Sec.~\ref{Sec:KRONIC} and~\ref{Sec:Regression}).}

\subsection{Dynamic mode decomposition}\label{Sec:Background:DMD} 
Dynamic mode decomposition (DMD) is a simple numerical algorithm that approximates the Koopman operator with a best-fit linear model that advances measurements from one time step to the next~\cite{Rowley2009jfm,Schmid2010jfm,Tu2014jcd,Kutz2016book}.  
DMD was originally introduced in the fluid dynamics community~\cite{Schmid2010jfm,Schmid:2011,Grilli:2012,Bagheri2013jfm,Mezic2013arfm,Brunton2015amr} to decompose large data sets into dominant spatial-temporal coherent structures, and a connection with Koopman theory was soon established~\cite{Rowley2009jfm}.   
In addition to fluid dynamics, DMD has been widely applied to a range of problems {in neuroscience~\cite{brunton2016extracting}, robotics~\cite{Berger2014ieee}}, epidemiology~\cite{Proctor2015ih}, and video processing~\cite{Grosek2014arxiv,kutzRPCA1,Erichson2016jrtp}.  
The DMD algorithm has also been extended to include sparsity and compressed sensing~\cite{Tu2014ef,Jovanovic2014pof,Brunton2015jcd,Gueniat2015pof,Erichson2016jrtp}, actuation and control~\cite{Proctor2016siads}, multi-resolution analysis~\cite{Kutz2016siads}, de-noising~\cite{Dawson2016ef,Hemati2015arxiv}, and streaming variants~\cite{Hemati2014pof}.   

In DMD, the infinite-dimensional Koopman operator in \eqref{Eq:KoopmanDiscrete2} is approximated with a finite-dimensional matrix $\bA$ that advances the system state $\bx$:
\begin{equation}
\bx_{k+1} \approx \bA \bx_k.\label{Eq:DMDLinear}
\end{equation}

Given data from a nonlinear system in \eqref{Eq:DiscreteDynamics}, it is possible to stack snapshots into a matrix 
$\bX = [\bx_1\; \bx_2\; \ldots\; \bx_{m-1}]$
and a time-shifted matrix 
$\bX^{'} = [\bx_2\; \bx_3\; \ldots\; \bx_m]$.  
In terms of these data matrices, \eqref{Eq:DMDLinear} becomes
\begin{equation}
\bX^{'} \approx \bA \bX.
\end{equation}
Various DMD algorithms then compute the leading eigendecomposition of the best-fit linear operator $\bA$, given by
\begin{equation}
\bA  = \mathrm{arg}\,\min\limits_{\tilde{\bA}} \vert\vert \bX^{'} - \tilde{\bA}\bX \vert\vert_{F}, 
\end{equation}
where $\|\cdot\|_F$ is the Frobenius norm.  
The best-fit $\bA$ is given by $\bA  = \bX^{'}\bX^{\dagger}$, where $\dagger$ is the pseudo-inverse, which is computed via singular value decomposition.

DMD has proven to be an extremely useful technique for the analysis of high-dimensional dynamical systems data.  
However, DMD is based on linear measurements of the system, which do not typically span a Koopman-invariant subspace of a general nonlinear system.  
For example, a linear DMD model may perfectly capture the periodic attractor dynamics of a system on a limit cycle, but will fail to capture the nonlinear transients if the system is perturbed off the attractor.  
DMD has since been augmented with nonlinear measurements to enrich the model in EDMD~\cite{Williams2015jnls,Williams2015jcd,klus2015numerical} and VAC~\cite{noe2013variational,nuske2016variational}.  
EDMD models have been used with success for estimation~\cite{Surana2016cdc,Surana2016nolcos} and model predictive control~\cite{Korda2016arxiv}.  
However, EDMD models are based on a large set of nonlinear measurements of the state, and there is no guarantee that these measurements form a Koopman invariant subspace.  
In fact, EDMD measurement subspaces will generally \emph{not} be closed~\cite{Brunton2016plosone,Kutz2016book}.    
For example, there is no finite-dimensional Koopman invariant subspace that includes the state of the system $\bx$ for any dynamical system that has multiple attractors (e.g., fixed points, periodic orbits, etc.), since the resulting finite-dimensional linear model cannot be topologically conjugate to the original dynamics.  
This is closely related to the representation of the Koopman operator in a polynomial basis, similar to Carleman linearization~\cite{banks1992infinite,kowalski1991nonlinear,svoronos1994discretization}.  
Thus, EDMD as well as other models are often plagued with spurious eigenfunctions that do not behave linearly as predicted by the associated eigenvalue.  
Fortunately, although these models may have corrupted eigenvalues and eigenfunctions, eigenfunctions corresponding to lightly damped eigenvalues may be faithfully extracted.

\subsection{Optimal control}\label{Sec:Background:OptimalControl}

The overarching goal of Koopman control is to reformulate strongly nonlinear dynamics in a linear framework to enable the use of powerful optimal and robust control techniques available for linear systems~\cite{sp:book,dp:book,stengel2012book}.  
Here, we summarize key results in optimal control theory that will be used for control in Koopman eigenfunction coordinates.  
We consider the nonlinear system~\eqref{Eq:Dynamics} affected by an external input
\begin{equation}\label{Eqn:OptimalControl:NonlinearSystemWithControl}
\frac{\mathrm d}{\mathrm dt} \bx(t) = {\boldsymbol{\it f}}(\bx,\bu),\quad \bx(0) = \bx_0
\end{equation}
with multi-channel control input $\bu\in\mathbb{R}^{q}$ and continuously differentiable dynamics
$ {\boldsymbol{\it f}}(\bx,\bu):\mathbb{R}^n\times \mathbb{R}^q\rightarrow\mathbb{R}^n$.
Without loss of generality, the origin is an equilibrium: ${\boldsymbol{\it f}}({\bf 0},{\bf 0}) = {\bf 0}$.

Infinite-horizon optimal control minimizes the following quadratic cost functional \begin{equation}\label{Eqn:OptimalControl:NL:QuadraticCostFunction}
J(\bx, \bu) = \frac{1}{2} \int\limits_0^{\infty}\, \underbrace{\bx^T(t) \bQ \bx(t) + \bu^T(t)\bR \bu(t)}_{L(\bx,\bu)} \,  d t
\end{equation}
with state and input weight matrices $\bQ\in\mathbb{R}^{n\times n}$ and $\bR\in\mathbb{R}^{q\times q}$.
Both matrices are symmetric and fulfill $\bQ > 0 $ and $\bR \geq 0$. 
A full-state feedback control law
\begin{equation}\label{Eqn:OptimalControl:NL:ControlLaw}
\bu(\bx) = - \bC (\bx)\bx = \bc(\bx),\quad \bc({\bf 0}) = {\bf 0}
\end{equation}
with gain $\bC:\mathbb{R}^n\rightarrow\mathbb{R}^{q\times n}$
is sought that minimizes the cost function~\eqref{Eqn:OptimalControl:NL:QuadraticCostFunction} 
subject to the state dynamics~\eqref{Eqn:OptimalControl:NonlinearSystemWithControl} to drive the system to the origin, i.e.\ $\lim\limits_{t\rightarrow\infty}\bx(t) = {\bf 0},\forall \,\bx$.

%% file: KRONIC.tex
We now propose a general control architecture in Koopman eigenfunction coordinates, referred to as \emph{Koopman Reduced Order Nonlinear Identification and Control} (KRONIC)(see Fig.~\ref{Fig:Schematic:KRONIC}).
\RefereeOne{
This eigenfunction perspective relies on a model constructed from validated Koopman eigenfunctions, which is closed and linear by design.	
}	
First, we derive how the effect of actuation affects the dynamics of these eigenfunction coordinates.
Second, the optimal control problem is formulated in these coordinates and a corresponding feedback controller is then developed, yielding a possibly nonlinear control law in the original state variables. 
Control in eigenfunction coordinates is quite general, encompassing the stabilization of fixed points and periodic orbits, e.g.\ via the Hamiltonian eigenfunction, or the manipulation of more general spatial-temporal coherent structures given by level sets of other eigenfunctions.

\subsection{Control--affine systems}
We first examine how adding control to a dynamical system \eqref{Eq:Dynamics} affects a single Koopman eigenfunction.  
This formulation then readily generalizes to multiple eigenfunctions. 
Consider a control-affine system 
\begin{equation}\label{Eqn:KRONIC:AffineControlSystem0}
\frac{\mathrm d}{\mathrm dt}\bx(t) = {\boldsymbol{\it f}}(\bx) + \sum_{i=1}^{q}{\bf b}_i(\bx)u_i,
\end{equation}
with a multi-channel input $\bu\in\mathbb{R}^{q}$,  continuously differentiable dynamics
$ {\boldsymbol{\it f}}(\bx):\mathbb{R}^n\rightarrow\mathbb{R}^n$ associated with the unforced dynamics, and each ${\bf b}_i(\bx)$ is a vector field acting on the state space. 

Starting with the Koopman operator associated with the uncontrolled, autonomous system (see Sec.~\ref{Sec:Background:Koopman}), we examine how the control terms in~\eqref{Eqn:KRONIC:AffineControlSystem0} affect the dynamics of its eigenfunctions. 
By applying the chain rule, we obtain
\begin{subequations}\label{Eq:KoopmanControl}
	\begin{align}
		\frac{\mathrm d}{\mathrm dt}\varphi(\bx) &= \nabla\varphi(\bx)\cdot\left({\boldsymbol{\it f}}(\bx) + \sum_{i=1}^{q}{\bf b}_i(\bx)u_i\right)\\
		&= \lambda\varphi(\bx) + \nabla\varphi(\bx)\cdot\sum_{i=1}^{q}{\bf b}_i(\bx)u_i.
	\end{align}
\end{subequations}
This equation differs from Eq.~\eqref{Eq:KoopmanEfun} in the additional second term associated with the control terms.
The $\varphi(\bx)$ is a Koopman eigenfunction associated with the autonomous Koopman operator for the unforced dynamics ${\bf f}(\bx)$.
For instance, a Hilbert space of the Lebesque square-integrable functions may be considered as function space.
The control enters the dynamics of $\varphi$ via the additional term leading to a control-affine system, which is linear in $\varphi$ and possibly nonlinear in the control.
	
Without loss of generality, we assume in the examples presented in later sections a linear control term:
\begin{equation}\label{Eqn:KRONIC:AffineControlSystem}
	\frac{\mathrm d}{\mathrm dt}\bx(t) = {\boldsymbol{\it f}}(\bx) + \bB \bu,
\end{equation}
with control matrix $\bB\in\mathbb{R}^{n\times q}$, so that the dynamics of the eigenfunctions become
\begin{subequations}
\begin{align}
\frac{\mathrm d}{\mathrm dt}\varphi(\bx)
&= \lambda\varphi(\bx) + \nabla\varphi(\bx)\cdot \bB\bu.
\end{align}
\end{subequations}
%

\subsection{Nonaffine control systems}\label{Sec:KoopmanEigenfunctionPDE_Nonaffine}
More generally, we may be interested in the control of a nonlinear, non-affine system:
\begin{equation}\label{Fig:KRONIC:NonAffineControlSystem}
\frac{\mathrm d}{\mathrm dt}\bx(t) = {\boldsymbol{\it f}}(\bx,\bu),
\end{equation}
with continuously differentiable dynamics $ {\boldsymbol{\it f}}(\bx,\bu):\mathbb{R}^n\times \mathbb{R}^q\rightarrow\mathbb{R}^n$.
We may now consider an observable $g(\bx,\bu)$ as a function of the extended state space $\mathcal{M}\times\mathcal{U}$, where $\bx\in\mathcal{M}$ and $\bu\in\mathcal{U}$, and the non-autonomous Koopman operator acting on these observables. 
For discrete-time dynamics with discrete map ${\bf F}(\bx,\bu)$, the Koopman operator propagates scalar measurements according to $Kg(\bx_k,\bu_k) = g({\bf F}(\bx_k,\bu_k),\bu_{k+1})$ by assuming the Koopman operator acts on the extended state space in the same manner as the Koopman operator associated with the autonomous, unforced system.
This assumption has been previously~\cite{Korda2016arxiv} considered, where the extended space was defined as the product of the original state-space and the space of all control sequences.
Since the first appearance of this article, they have been further studies~\cite{Proctor2016arxiv, Arbabi2019cdc,Kaiser2020springer} examining the approximation of the Koopman operator for the non-affine control system.
In~\cite{Proctor2016arxiv}, it is in more detail discussed how the Koopman operator formulation can be modified based on the dynamics of $\bu$ itself, e.g. open-loop versus closed-loop control.
As associated function space, in which observables are defined on the extended state space, a Hilbert space of the Lebesque square-integrable functions or polynomial functions defined on a compact set may be considered.
If the dynamics on $\bu$ are governed by a specific feedback law of the form $\bu = \bC (\bx)$, that is only a function of the state $\bx$, then choices about the function space for the controlled system are equivalent to those applying to the autonomous system. 

For smooth dynamics~\eqref{Fig:KRONIC:NonAffineControlSystem}, the continuous-time dynamics are given by
\begin{equation}
\frac{d}{dt}g(\bx,\bu) = \mathcal{K}g(\bx,\bu).
\end{equation}
An associated Koopman eigenfunction $\varphi(\bx,\bu)$ satisfies
\begin{align}
\frac{\mathrm d}{\mathrm dt}\varphi(\bx,\bu) &= \lambda\varphi(\bx,\bu).
\end{align}
Applying the chain rule, we find that the dynamics of the Koopman eigenfunction depends on $\dot{\bu}$, which is generally arbitrary:
$\frac{\mathrm d}{\mathrm dt}\varphi(\bx,\bu) = \nabla_{\bx}\varphi(\bx,\bu)\cdot\mathbf{f}(\bx,\bu) + \nabla_{\bu}\varphi(\bx,\bu)\cdot\dot{\bu}$.  
Instead, we may specify that $\varphi(\bx,\bu)$ reduces to the eigenfunction $\varphi(\bx,\bar{\bu})$ of $\dot{\bx} = {\boldsymbol{\it f}}(\bx,\bar{\bu})$ for all locked $\bar{\bu}\in\mathcal{U}$, as in~\cite{Proctor2016arxiv}.  
In this case, the eigenfunction is \emph{parametrized} by the input $\bar{\bu}$
\begin{equation}
\nabla_{\bx}\varphi(\bx,\bar{\bu})\cdot{\boldsymbol{\it f}}(\bx,\bar{\bu}) = \lambda\varphi(\bx,\bar{\bu}),\quad\forall ~\bar{\bu}\in\mathcal{U}.
\end{equation}
These are eigenfunctions of the parametrized Koopman generator, which is autonomous for each locked $\bar{\bu}$ and can be defined on the commonly used function spaces as stated above.
\RefereeOne{
This perspective is also assumed within the gEDMD framework and its extension for control~\cite{Klus2020arxiv}.
The Koopman generator $\mathcal{K}$ is approximated as a finite-rank matrix $\bK$ parametrized by the discrete control input using EDMD.
The control problem is then solved by optimizing switching times among the finite set of discrete control inputs and associated models.}

If we augment the eigenfunction vector with the input $\bu$, we obtain
\begin{equation}\label{Eqn:KRONIC:GainScheduled}
\frac{\mathrm d}{\mathrm dt}\varphi(\bx,\bu) = \lambda\varphi(\bx,\bu) + \nabla_{\bu}\varphi(\bx,\bu)\cdot\dot{\bu},
\end{equation}
where we view $\dot{\bu}$ as the input to the Koopman linear system, and the $\nabla_{\bu}\varphi(\bx,\bu)$ matrix varies based on $\bx$ and $\bu$.    
Thus, we may enact a gain-scheduled control law.

Summarizing, 
if the original dynamics are nonlinear and control-affine, the dynamics in Koopman eigenfunction coordinates are control-affine and split into a linear part associated with the nonlinear unforced dynamics and a bilinear part associated with the control term.
If the original dynamics are nonlinear and non-affine in the control, the eigenfunction dynamics can become linear if the Koopman operator is defined on the extended state. Considering practical implications, it is also possible to modify these dynamics so that these consist of a linear part and a bilinear control term as in Eq.~\eqref{Eqn:KRONIC:GainScheduled}.
We also point out, that	while the control is generally nonlinear in the state $\bx$, it may become linear in the eigenfunction coordinates for special cases, such as $\nabla\varphi(\bx) = \mathrm{const.}$

\subsection{Formulation of the optimal control problem}\label{Sec:KRONIC:OptimalControl}
We now formulate the infinite-horizon, optimal control problem~\cite{stengel2012book} for a reduced set of Koopman eigenfunctions.
The control objective is a quadratic cost functional: 
\vspace{-0.1in}
\begin{equation}\label{Eqn:KRONIC:QuadraticCostFunction}
\vspace{-0.1in}
\hspace{-0.1in}J(\bvarphi, \bu) = 
\frac{1}{2} \int\limits_0^{\infty}\, \bvarphi^T(\bx(t)) \bQ_{\varphi} \bvarphi(\bx(t)) + \bu^T(t)\bR \bu(t) \, \mathrm d t,
\end{equation}
where $\boldsymbol{\varphi} =[\varphi_{1}\,\varphi_{2}\,\ldots\,\varphi_{r}]^T$ comprises $r$ eigenfunctions with $\varphi_{j}$ associated with eigenvalue $\lambda_j$.
For this cost function to be equivalent to the cost in the original state $\bx$, a modified weight matrix may be considered such that $\bvarphi^T \bQ_{\varphi} \bvarphi \approx \bx^T\bQ \bx$. 
This can only be achieved exactly if the state itself is an eigenfunction of the Koopman operator; 
however, it can be sufficient requiring this only for a subset of states that enter the cost function.
Alternatively, it is possible to estimate $\bx$ via the inverse mapping $\varphi^{-1}$. 
Eigenfunctions may generally not be invertible exactly.
Nevertheless, the inverse mapping may be approximated using multidimensional scaling as in~\cite{kawahara2016nips} or learned jointly with $\varphi$ itself using autoencoders~\cite{Lusch2018natcomm}.
More generally, 
the matrix $\bQ_{\varphi}$ allows one to weight particular eigenfunction directions, which are related to properties of the underlying system and coherent structures.
The selection of a specific set of eigenfunctions, in which a model is constructed and which are used to formulate the cost functional, is problem specific.
However, given a target state $\bx^{REF}$ in the original state space, the associated target value of the eigenfunctions may be directly determined by evaluating the eigenfunctions on the target state, i.e. $\bvarphi^{REF} := \bvarphi(\bx^{REF})$.

For the general case, 
it is possible to augment the state with the control input and include the derivative of the control as new input $\hat{\bu}:=\dot{\bu}$:
\begin{equation}
\frac{\mathrm d}{\mathrm dt}
\begin{bmatrix}
\bvarphi\\ \bu
\end{bmatrix}
=
\begin{bmatrix}
\bLambda & \bB_{\varphi}\\ {\bf 0} & {\bf 0}
\end{bmatrix}
\begin{bmatrix}
\bvarphi\\ \bu
\end{bmatrix}
+
\begin{bmatrix}
{\bf 0}\\ {\bf I}_{q}
\end{bmatrix}
\hat{\bu}
\end{equation}
with $q\times q$ identity matrix ${\bf I}_q$.  This may be interpreted as integral control.  
The cost functional is then given by
\begin{equation}
J = \frac{1}{2} \int\limits_{0}^{\infty}\, 
\begin{bmatrix}
\bvarphi^T & \bu^T
\end{bmatrix}
\begin{bmatrix}
\bQ_\varphi & {\bf 0}\\ {\bf 0} & {\bf R}
\end{bmatrix}
\begin{bmatrix}
\bvarphi\\ \bu
\end{bmatrix}
+
\hat{\bu}^T \hat{\bf R} \hat{\bu} \, \mathrm d t.
\end{equation}
with some restrictions on $\hat{\bR}$. 
Modifying the system structure, by moving the nonlinearity in the control term into the state dynamics, 
improves the tractability of the problem~\cite{beeler2004report}.

In the following, we will focus on multiple-input, control-affine systems~\eqref{Eqn:KRONIC:AffineControlSystem}, for which the dynamics in intrinsic coordinates becomes
\begin{equation}\label{Eqn:KRONIC:MIMO}
\frac{\mathrm d}{\mathrm dt}\boldsymbol{\varphi}(\bx) = \bLambda\boldsymbol{\varphi}(\bx) + \nabla_{\bx}\boldsymbol{\varphi}(\bx)\cdot\bB\bu
\end{equation}
with $\bLambda = \mathrm{diag}({\lambda_1},\,\ldots,\,{\lambda_r})$.  
Depending on the structure of $\bvarphi(\bx)$ and $\bB$, the actuation matrix $\bB_{\varphi} = \nabla_{\bx}\boldsymbol{\varphi}(\bx)\cdot\bB$ may be a function of $\bx$. 
A state-dependent control term may be interpreted as a gain-scheduled control. 
A feedback controller is now sought in the Koopman representation of the form
\begin{equation}\label{Eqn:KRONIC:FeedbackControlLaw}
\bu =-\bc_{\varphi}[\bvarphi(\bx)] =  -\bC_{\varphi}(\bx) \bvarphi(\bx).
\end{equation}
We may also consider reference tracking,
$\bu = -\bC_{\varphi}(\bx) \left[\bvarphi(\bx) - \bvarphi^{REF}\right]$, with a  modified cost functional~\eqref{Eqn:KRONIC:QuadraticCostFunction}.
The objective is then to determine the gain function $\bC_{\varphi}$ by minimizing $J(\bvarphi,\bu)$ and the resulting control $\bu$ may be directly applied to the original system~\eqref{Eqn:KRONIC:AffineControlSystem}, for which the control problem may be suboptimal.
Generally, Koopman eigenfunction control can be used in two ways depending on which means of control one has access to: 
(1) internally driven swimmers or particles, which state is given by $\bx$ and which dynamics are subject to an external field, such as a fluid flow or magnetic field; or 
(2) driving an external field, represented in terms of these eigenfunctions, in which these swimmers or particles drift.

\subsection{Solving the optimal control problem}
The model~\eqref{Eqn:KRONIC:MIMO} may be combined with any model-based control strategy. In the following, we discuss ways to solve the optimal control problem formulated above using standard techniques.
	
In the simplest case, for which Eq.~\eqref{Eqn:KRONIC:MIMO} becomes fully linear,
the optimal gain matrix $\bC_{\varphi}$ can be determined by solving the associated algebraic Riccati equation leading to a linear quadratic regulator (LQR) formulated in Koopman eigenfunctions. 
For the more general case, where the control term becomes nonlinear, other techniques are required. 
A common extension of LQR for nonlinear systems considers state-dependent state and control matrices as outlined in Sec.~\ref{Sec:Control:SDRE} and solves the state-dependent Riccati equation~\eqref{Eqn:OptimalControl:SDRE}. Here, the state dynamics would be constant and linear, i.e. $\bA = \bLambda$, and only the control matrix $\bB_{\varphi}:= \nabla_{\bx}\varphi\cdot \bB$ depends on the state $\bx$:
\begin{equation}\label{Eqn:KRONIC:SDRE}
\bQ + \bP\bLambda + \bLambda^T\bP - \bP\bB_{\varphi}(\bx)\bR^{-1}\bB_{\varphi}^T(\bx)\bP = {\bf 0}.
\end{equation}
We note that parametrizations for the dynamics~\eqref{Eqn:OptimalControl:FacorizedLinearSystem} in $\bx$ are unique for scalar systems, but generally nonunique for multivariable systems~\cite{Cloutier1996proc}.
In contrast, the state-dependent dynamics in eigenfunction coordinates~\eqref{Eqn:KRONIC:MIMO} are unique with respect to $\bx$ and $\bLambda$ is constant. These are not a result from the factorization.
The (non)uniqueness of the factorization is generally related to global optimal control and global optimal stability. However, further studies are required to connect these properties for solutions in eigenfunction coordinates to the state dynamics in $\bx$.

\RefereeOne{
We examine the effect of an error $\varepsilon\psi(\bx)$ in the representation of a Koopman eigenfunction, $\hat\varphi(\bx) := \varphi(\bx) + \varepsilon\psi(\bx)$,  on its closed-loop dynamics based on~\eqref{Eqn:KRONIC:SDRE} and provide an upper bound for the error (for details see Appendix~\ref{Sec:ControlError}). 
We assume control-affine dynamics~\eqref{Eqn:KRONIC:AffineControlSystem} of the underlying system, access to full-state measurements $\bx$, and control vector fields are known. Further, we reformulate $\bB(\bx)\bu:=\sum_{i=1}^q \bb_i(\bx)u_i$ for simplicity.
The closed-loop dynamics of $\varphi(\bx)$ are then given by
\begin{equation}
\dot{\varphi}(\bx) = -\sqrt{\lambda^2+Q\bC(\bx)\bR^{-1}\bC^T(\bx)}\varphi(\bx) = -\sqrt{\mu}\varphi(\bx),
\end{equation}	
where $\mu := \lambda^2+Q\bC(\bx)\bR^{-1}\bC^T(\bx)$.
The upper bound for the error in $\mu$ due to the misrepresentation of $\varphi(\bx)$ due to $\varepsilon\psi(\bx)$is
\begin{equation}
\vert\mu-\hat{\mu} \vert 
\leq \Big\lvert-\varepsilon \left(Q\bC\bR^{-1}\bD^T+Q\bD\bR^{-1}\bC^T+2\lambda\frac{\nabla\psi\cdot{\bf f}}{\varphi(\bx)}\right)-\varepsilon^2\left(Q\bD\bR^{-1}\bD^T+\frac{(\nabla\psi\cdot{\bf f})^2}{\varphi^2(\bx)}\right)\Big\rvert.
\end{equation}
where
$\bC(\bx):=\nabla\varphi(\bx)\cdot\bB(\bx)$ and  $\bD(\bx):=\nabla\psi(\bx)\cdot\bB(\bx)$ are the control terms associated with $\varphi(\bx)$ and $\psi(\bx)$, respectively.
For small $\varphi(\bx)$ the contribution of $\varepsilon\psi(\bx)$ becomes important/may be dominant.
}

\RefereeOne{
More generally, Koopman control in eigenfunction coordinates may be combined with any model-based control approach.}
Under certain conditions, it may also be possible to feedback linearize the dynamics~\cite{khalil1996book}.
The data-driven identification of Koopman eigenfunctions can be challenging and they may only be approximated accurately in a certain domain. 
Further, dynamics may also drift away from the situations captured in the training data due to external disturbances.
Especially in these cases it is advantageous to couple the resulting model with a receding horizon estimator and controller to adapt quickly to changing conditions. 
In particular, model predictive control has gained increasing popularity over the last decade 
due to its success in a wide range of applications and its ability to incorporate customized cost functions and constraints~\cite{camacho2013model,allgower2004nonlinear,eren2017jgcd}.

%% file: Regression.tex
The overarching goal of this work is to control nonlinear systems in intrinsic Koopman eigenfunction coordinates.
\RefereeOne{
A leading method for approximating the Koopman operator is EDMD relying on a set of basis or dictionary functions.
Here, we aim to address the well-known closure issue, i.e.\ the model may not be closed in the set of basis functions and can therefore have spurious eigenfunctions, by assuming an eigenfunction perspective.
Approximating these eigenfunctions from data is an ongoing challenge.  
Building on EDMD, we propose a strategy to directly identify dominant Koopman eigenfunctions associated with lightly damped eigenvalues, that can then be utilized to construct low-dimensional, closed models.
Further, we show that it is possible to sparsify the state-transition matrix $\bK$ of EDMD to improve the representation of eigenfunctions within the selected basis. }

\RefereeALL{
It has long been recognized that numerical regularization via the truncation of the SVD in the computation of the pseudoinverse, as it is done for DMD/EDMD, is crucial to reduce noise corruption.
However, this does not produce a sparser approximation of $\bK$, which is crucial for systems that are inherently sparse in their representation.
An additional $L1$ regularization is here advantageous as it is able to improve the predictive power of the model and eigenfunctions become more accurately represented in the chosen function library preventing overfitting.}
\RefereeOne{
Interestingly, the infinitesimal generator of the Koopman operator may be sparse even when the Koopman operator itself is not~\cite{Klus2020arxiv}. However, even in these cases, the approach may benefit from a sparsity constraint to counteract spurious non-zero entries arising from noise and numerical approximation. }
\RefereeALL{
More generally, non-compactness or continuous spectra of the Koopman operator can pose issues for the numerical analysis requiring some form of regularization (see~\cite{Giannakis2019} for a detailed discussion for the Koopman operator and~\cite{Froyland2013} regarding the Perron-Frobenius operator).}


In the following we formulate a framework using sparsity-promoting techniques to identify Koopman eigenfunctions directly building on the partial differential equation (PDE) governing the evolution of an eigenfunction.
Applying the chain rule to \eqref{Eq:KoopmanEfun} yields
\begin{equation}
\frac{\mathrm d}{\mathrm dt}\varphi(\bx) = \nabla\varphi(\bx)\cdot \dot{\bx} = \nabla\varphi(\bx)\cdot \mathbf{f}(\bx).\label{Eq:KoopmanDynamics}
\end{equation}
Combined with \eqref{Eq:KoopmanEfun}, this results in a linear PDE for the eigenfunction $\varphi(\bx)$:
\begin{equation}
{\nabla\varphi(\bx)\cdot\mathbf{f}(\bx)= \lambda\varphi(\bx)}.\label{Eq:KoopmanPDE}
\end{equation}
This formulation assumes continuous and differentiable dynamics and that the eigenfunctions are smooth~\cite{Mezic2017book}.

Spectral properties of the Koopman operator have been shown to relate to intrinsic time scales, geometrical properties, and the long-term behavior of dynamical systems~\cite{Mezic2005nd,Mezic2004physicad,mauroy2013cdc,Kutz2016book}.
It has been shown~\cite{Mezic2017arxiv}, that the evolution of observables can be described by a linear expansion in Koopman eigenfunctions for systems which consist only of the point spectrum, i.e. smooth dynamical systems exhibiting, e.g., hyperbolic fixed points, limit cycles and tori.
If systems with a mixed spectrum are considered, it may be possible to restrict the following analysis to the point spectrum as in~\cite{Das2019jsp}. 

\subsection{Data-driven discovery of continuous-time eigenfunctions}\label{Sec:Regression:CT}
\begin{figure}[tb]
	\begin{center}
		\includegraphics[width=\textwidth]{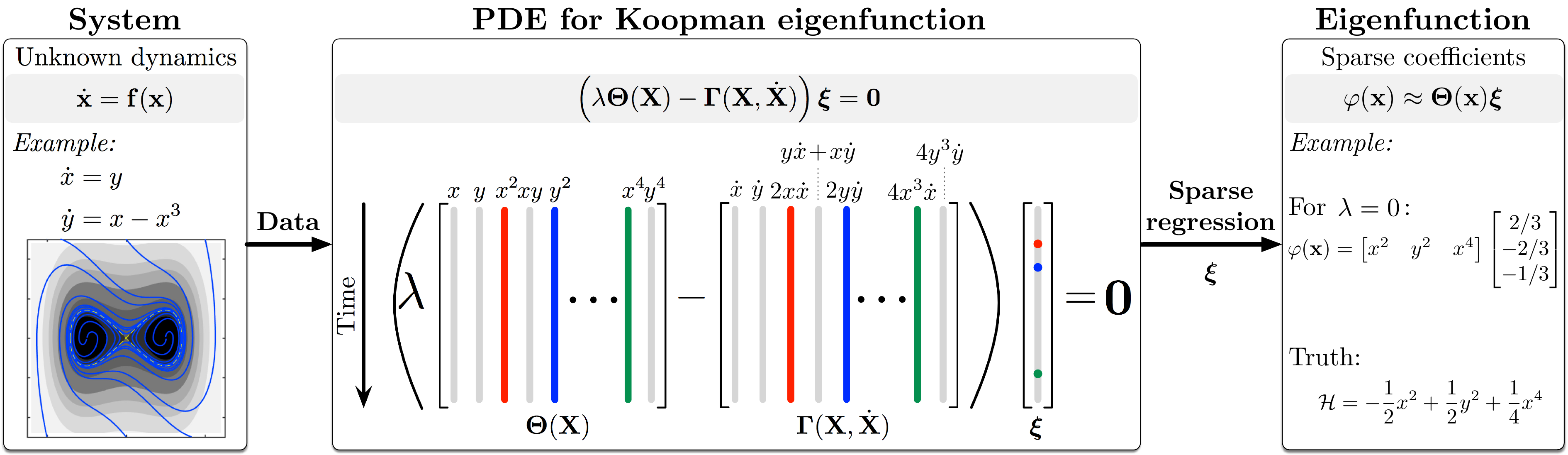}
		\caption{Identification of Koopman eigenfunctions from data using implicit-SINDy.}
			\vspace{-.1in}
		\label{Fig:ImplicitSINDy}
	\end{center}
\end{figure}

Sparse identification of nonlinear dynamics (SINDy)~\cite{Brunton2016pnas} is used to identify Koopman eigenfunctions for a particular value of $\lambda$.  
This formalism assumes that the system has a point or mixed spectrum, for which eigenfunctions with distinct eigenvalues exist.
A schematic is displayed in Fig.~\ref{Fig:ImplicitSINDy}.

First, we build a library of candidate functions:
\begin{equation}
\bTheta(\bx) = \begin{bmatrix} \theta_1(\bx) &\theta_2(\bx) & \cdots & \theta_p(\bx)\end{bmatrix}.
\end{equation}
We choose $\bTheta$ large enough so that the Koopman eigenfunction may be well approximated in this library:
\begin{equation}
\varphi(\bx) \approx \sum_{k=1}^p\theta_k(\bx)\xi_k = \bTheta(\bx)\bxi.
\end{equation}
Given data $\bX = [\bx_1\, \bx_2\, \cdots\, \bx_m]$, the time derivative  $\dot{\bX}=[\dot{\bx}_1\, \dot{\bx}_2\, \cdots\, \dot{\bx}_m]$ can be approximated numerically from $\bx(t)$ if not measured directly~\cite{Brunton2016pnas}. 
The total variation derivative~\cite{Chartrand2011isrnam} is recommended for noise-corrupted measurements.
It is then possible to build a data matrix $\bTheta(\bX)$:
%
\begin{equation}\label{Eq:Theta}
\bTheta(\bX) = \begin{bmatrix} \theta_1(\bX^T) & \theta_2(\bX^T) & \cdots & \theta_p(\bX^T)\end{bmatrix}.
\end{equation}
Moreover, we can define a library of directional derivatives, representing the possible terms in $\nabla\varphi(\bx)\cdot\mathbf{f}(\bx)$ from~\eqref{Eq:KoopmanPDE}:  $\bGamma(\bx,\dot{\bx})=[\nabla\theta_1(\bx)\cdot\dot{\bx}\, \nabla\theta_2(\bx)\cdot\dot{\bx}\, \cdots\, \nabla\theta_p(\bx)\cdot\dot{\bx} ]$.  It is then possible to construct $\bGamma$ from data:
\begin{align}\label{Eq:Gamma}
\bGamma(\bX,\dot{\bX}) =
\begin{bmatrix} \nabla\theta_1(\bX^T)\cdot\dot{\bX} & \nabla\theta_2(\bX^T)\cdot\dot{\bX} & \cdots & \nabla\theta_p(\bX^T)\cdot\dot{\bX}\end{bmatrix}.
\end{align}

For a specific eigenvalue $\lambda$, the Koopman PDE in \eqref{Eq:KoopmanPDE} may be evaluated on data, yielding:
\begin{equation}
\left(\lambda\bTheta(\bX) - \bGamma(\bX,\dot{\bX})\right)\bxi = \mathbf{0}.\label{Eq:SparseKoopman}
\end{equation}

The formulation in \eqref{Eq:SparseKoopman} is implicit, so that $\bxi$ will be in the null-space of the matrix ${\lambda\bTheta(\bX)-\bGamma(\bX,\dot{\bX})}$.  
The right null-space of~\eqref{Eq:SparseKoopman} for a given $\lambda$ is spanned by the right singular vectors of $\lambda\bTheta(\bX)-\bGamma(\bX,\dot{\bX}) = \bU \bSigma\bV^{*}$ (i.e., columns of $\bV$) corresponding to zero-valued singular values.  
It is possible to identify the few active terms in an eigenfunction by finding the sparsest vector in the null-space~\cite{Qu2014}, which is used, e.g., in the implicit-SINDy algorithm~\cite{Mangan2016ieee}. This is a nonconvex approach based on alternating directions (adm) with linear scaling~\cite{Qu2014}, which can be adapted to our problem.
In this formulation, the eigenvalues $\lambda$ are not known \emph{a priori}, and must be learned online along with the approximate eigenfunction.  
	\begin{algorithm}
		\caption{Direct discovery of eigenfunctions based on the implicit formulation~\eqref{Eq:SparseKoopman}.}
		\label{Alg:Implicit}
		\begin{algorithmic}[1]
			\STATE{Initialize: ${\boldsymbol\Lambda}^{guess} = \mathrm{eigs}(\bTheta(\bX)^{\dagger}\bGamma(\bX,\dot{\bX}))$\}}
			\FOR{each eigenvalue $\lambda_i$ in ${\boldsymbol \Lambda}^{guess}$}
			\STATE{Initialize: $\lambda = \lambda_i^{guess}$}
			\WHILE{$\lambda$ not converged}
			\STATE{${\bf M} := \bGamma(\bX,\dot{\bX})-\lambda \bTheta(\bX)$}
			\STATE{${\bf N} := \mathrm{null}({\bf M})$}
			\FOR{each row $l$ in ${\bf N}$}
			\STATE{${\bf q}_0 := ({\bf N}(l,:)/\mathrm{norm}({\bf N}(l,:)))^T$}
			\STATE{${\bf q}_{mtx}(:,l) := \mathrm{adm}({\bf N},{\bf q}_0,\alpha,MaxIter,tolerance);$}
			\STATE{${\boldsymbol \Xi}(:,l) = {\bf N} {\bf q}_mtx(:,l)$}
			\ENDFOR
			\FOR{each column $l$ in ${\boldsymbol \Xi}$}
			\STATE{${\boldsymbol \Xi}(:,l):=\mathrm{soft\_thresholding}({\boldsymbol \Xi}(:,l) ,\alpha)$}
			\STATE{${\boldsymbol \Xi}(:,l):= \mathrm{normalize}({\boldsymbol \Xi}(:,l))$}
			\STATE{$Error(l) := \mathrm{norm}(\bTheta(\bX^{test}) {\boldsymbol \Xi}(:,l) - \mathrm{e}^{\lambda {\bf t}^{test}}\bTheta(\bx^{test}_0){\boldsymbol \Xi}(:,l))$}
			\ENDFOR
			\STATE{$best:= \mathrm{min}(Error)$}
			\STATE{${\boldsymbol{\xi}} := {\boldsymbol \Xi}(:,best)$}
			\STATE{Update $\lambda :=  ({\boldsymbol{\xi}}^T \bTheta(\bX) ^{\dagger} \bGamma(\bX,\dot{\bX}){\boldsymbol{\xi}} ) / ({\boldsymbol{\xi}}^T{\boldsymbol{\xi}})  $}
			\ENDWHILE
			\STATE{${\bf T}(:,i) =  {\boldsymbol{\xi}}$}
			\STATE{${\boldsymbol\Lambda}(i,i) = \lambda$}
			\ENDFOR
			\RETURN ${\boldsymbol\Lambda},{\bf T}$
		\end{algorithmic}
\end{algorithm}
In Alg.~\ref{Alg:Implicit}, we propose an implicit formulation, which starts with an initial guess of the eigenvalues given by the least-squares solution, and subsequently alternates between an searching for the sparsest vector in the null-space~\cite{Qu2014} and updating of the eigenvalue.
As the approach in~\cite{Qu2014} depends on the initial condition, we evaluate all initial conditions given by each row in the nullspace ${\bf N}:=\mathrm{null}(\bGamma(\bX,\dot{\bX})-\lambda \bTheta(\bX))$.
An additional soft-thresholding with parameter $\alpha$ and validation on a test dataset is applied to select the best eigenvector for each initial condition. 
While the approach has been observed to converge to accurate eigenvalues in clean data, the found eigenvector in the solution set may not be unique.
This approach solves for each eigenpair separately; however, it may also be possible to extend it to solve for eigenpairs jointly which is part of ongoing research.
From a practical standpoint, 
data in $\bX$ does not need to be sampled from full trajectories,
but can be obtained using more sophisticated strategies such as latin hypercube sampling or sampling from a distribution over the phase space.
It may also be possible to directly identify a recursion relationship to obtain a power series expansion as shown in Appendix~\ref{Sec:AnalyticalKoopmanEfun} and in~\cite{Mauroy2016ieee}.  

Koopman eigenfuntions and eigenvalues can also be determined as the solution to the eigenvalue problem $\bxi _{\alpha}\bK = \lambda_{\alpha}\bxi _{\alpha}$, where 
$\bK = \bTheta^{\dagger} \bGamma$ is obtained via least-squares (LS) regression.
While many eigenfunctions are spurious, i.e.\ these eigenfunctions do not behave linearly as predicted by the corresponding eigenvalue, those corresponding to lightly damped eigenvalues can be well approximated~\cite{Brunton2016plosone,Kutz2016book}, and a reduced-order Koopman model may developed on these coordinates (see also Sec.~\ref{Sec:Motivation}). 
The accuracy of eigenfunctions of $\bK$ can be improved, by improving the recovery of $\bK$ itself, which 
can be achieved by sparsifying $\bK$ directly. The minimization problem,
\begin{equation}\label{Eqn:Regression:SparseEDMD}
\min\limits_{{\bf k}_i}\, || \nabla\theta_i(\bX^T)\cdot\dot{\bX} -  \bTheta(\bX){\bf k}_i^{T}||^2_2 + 
\rho({\bf k}_i), \quad\forall {\bf k}_i, i=1,\ldots,p,
\end{equation}
where $\rho$ is a regularization term that promotes the sparsity of ${\bf k}_i$, i.e. the number of non-zero coefficients, 
is solved separately for each row in $\bK = [{\bf k}_1^T,\ldots,{\bf k}_p^T]$.
For instance, an $L_1$ constraint on the coefficients in ${\bf k}_i$ may be chosen and ${\bf k}_i$ may equivalently be determined as in SINDy~\cite{Brunton2016pnas}.
In general, the problem~\eqref{Eqn:Regression:SparseEDMD} can be solved using standard techniques, such as LASSO~\cite{SantosaLASSO,TibshiraniLASSO}, Least Angle Regression~\cite{BradleyLARS}, or an iterated least-squares thresholding method~\cite{Brunton2016pnas}.

\RefereeOne{
This formulation is closely related to the gEDMD framework~\cite{Klus2020arxiv}, a generalization of the EDMD method to approximate the infinitesimal generator of the Koopman operator. The gEDMD least-squares formulation is solved here row-wise with an additional sparsity constraint.}
There are also similarities to sparsity-promoting DMD~\cite{Jovanovic2014pof} and variants~\cite{Hua2017nd,Graff2019aiaa}, which aims to reconstruct a signal through a sparse combination of modes/eigenfunctions, which have been computed from a least-squares solution.
In contrast to these works, we argue that the sparse representation of eigenfunctions themselves should be promoted.
Alternatively, it is also possible to use the dominant terms in accurately identified eigenfunctions, generally associated with lightly damped eigenvalues,  as guidance to select observables to regress on.
Since the first appearance of this article, further promising methods have been proposed to identify Koopman eigenfunctions directly~\cite{Korda2020tac,Haseli2019arxiv,Pan2020arxiv}.

\subsection{\hspace{-0in}Data-driven discovery of discrete-time eigenfunctions}\label{Sec:Regression:DT}
In discrete-time, an eigenfunction evaluated at a number of data points in $\bX$ will satisfy:
\begin{equation}
\begin{bmatrix}
\lambda\varphi(\bx_1) &
\ldots &
\lambda\varphi(\bx_{m})
\end{bmatrix}^T = 
\begin{bmatrix}
\varphi(\bx_2) &
\ldots &
\varphi(\bx_{m+1})
\end{bmatrix}^T.
\end{equation}
Again, searching for such an eigenfunction $\varphi(\bx)$ in a library $\bTheta(\bx)$ yields the matrix system:
\begin{equation}\label{Eq:DiscreteEDMD}
\left(\lambda\bTheta(\bX)-\bTheta(\bX')\right)\bxi = \mathbf{0},
\end{equation}
where $\bX'= \begin{bmatrix}\bx_2 & \bx_3 & \cdots & \bx_{m+1}\end{bmatrix}$ is a time-shifted matrix. 
This formalism directly identifies the functional representation of a Koopman eigenfunction with eigenvalue $\lambda$.

If we seek the best \emph{least-squares} fit to~\eqref{Eq:DiscreteEDMD}, this reduces to the extended DMD~\cite{Williams2015jnls} formulation:
\begin{equation}\label{Eqn:KoopmanEfunLSR}
\lambda\bxi =  \bTheta^{\dagger}\bTheta' \bxi.
\end{equation}
Again, it is necessary to confirm that predicted eigenfunctions actually behave linearly on trajectories.  

\subsection{Koopman model reduction and validation}
\RefereeOne{
For any data-driven modeling framework cross-validation is critical to ensure predictive and generalization capabilities.
Approximations of the Koopman operator allow a systematic evaluation based on its linearity property.
We can learn a finite-rank approximation using any EDMD-like method or learn eigenfunctions directly.
In either case, it is critical to validate that any candidate eigenfunction $\varphi(\bx(t))$ actually behaves linearly on trajectories $\bx(t)$ as the eigenvalue $\lambda$ predicts.
The error of an eigenfunction can be defined as
\begin{equation}
E = ||\varphi(\bx(t)) - \mathrm{e}^{\lambda t}\varphi(\bx(0))||^2
\end{equation}  
evaluated on a test trajectory $\bx(t)$ and identified eigenfunctions can be ranked according to the error $E$. All eigenfunctions with error below a threshold may then be used to construct a reduced-order model:
\begin{equation}
\boldsymbol\varphi(\bx) = \boldsymbol\Lambda\boldsymbol{\varphi} + \nabla \boldsymbol\varphi(\bx)\cdot \bB\bu
\end{equation}
analogous to Eq.~\eqref{Eqn:KRONIC:MIMO}.
This model is closed and behaves linearly by design and we demonstrate in the following its increased predictive power.
}

%% file: SlowManifold.tex
\begin{figure}[!htb]
\centering
\includegraphics[width=\textwidth]{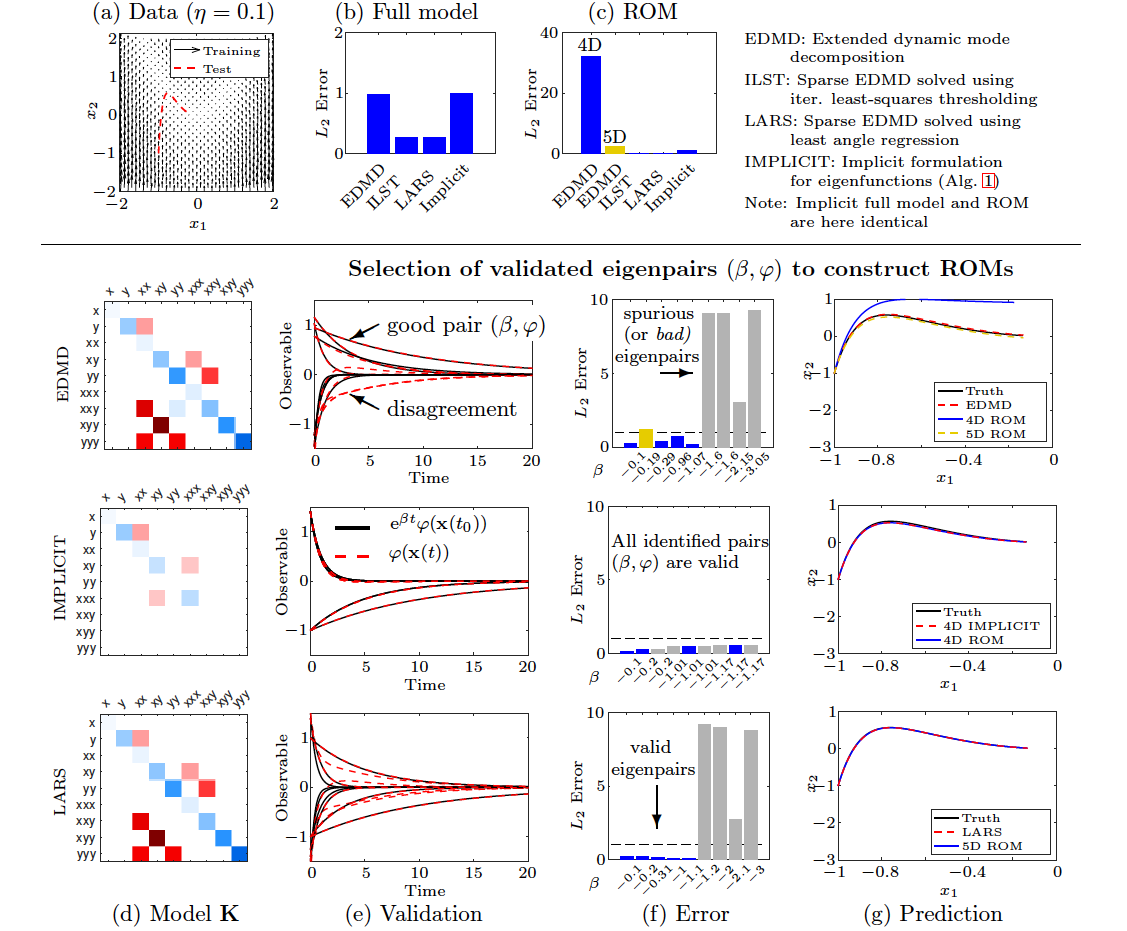}
	\caption{
	\RefereeOne{	
		Model comparison for the slow manifold system with noise level $\eta=0.1$:
		(a) training and test data,
		(b) state prediction error on test data for each identified model,
		(c) state prediction error on test data for each identified reduced-order model (ROM),
		(d) transition matrix on set of measurement functions, 
		(e) eigenpairs $(\beta,\varphi(\bx))$ validated on test trajectory (marked blue if below threshold (dashed line)), 
		(f) prediction error for each eigenpair $(\beta,\varphi(\bx))$ 
			(error between trajectory as predicted by eigenvalue and by evaluating eigenfunction on test trajectory), 
		(g) state prediction using model $\bK$ and ROM constructed from validated (or $\emph{good}$) eigenfunctions (marked in blue/yellow in (f)).
		For details see text.
	}
}
	\label{Fig:SLOWMANIFOLD:eta01}
\end{figure}

\subsection{Identification of eigenfunctions and models}
We consider a system with quadratic nonlinearity that gives rise to a slow manifold~\cite{Tu2014jcd}:
\begin{align}\label{Eqn:Motivation:SlowManifold1}
	\frac{\mathrm d}{\mathrm dt}\begin{bmatrix} x_1 \\ x_2 \end{bmatrix} = \begin{bmatrix} \mu x_1 \\ \lambda (x_2 - x_1^2)\end{bmatrix}
\end{align}
with $\mu=-0.1$ and $\lambda=-1$.
In the following, we apply the methods outlined in Sec.~\ref{Sec:Regression} to the dynamical system~\eqref{Eqn:Motivation:SlowManifold1}.
For all models, a library of the first $9$ monomials (up to the third degree)  is considered.
Noise corruption in measurement data can be particularly problematic. 
Here, we examine the recovery of eigenfunctions and prediction performance for different noise magnitudes $\eta=0.1$ and $\eta=0.9$.
Further, reduced-order models are constructed based on eigenpairs $(\beta,\varphi(\bx))$ that are deemed accurate, i.e. have a small $L_2$ error when compared with the prediction using the associated eigenvalue, here denoted by $\beta$.
The threshold for the selection of eigenfunctions is $1$ and $3$ for noise magnitudes $\eta=0.1$ and $\eta=0.9$, respectively.
The $L_1$ regularized problem is solved using the iterative least-squares thresholding algorithm (ITLS) as in SINDy~\cite{Brunton2016pnas} and least angle regression (LARS)~\cite{BradleyLARS}; however, detailed results are only shown for LARS as these have equivalent performance.
LARS has the advantage that the number of iterations scales with the number of candidate functions; only $10$ iterations are required here.
The sparse solution is then selected when there is minimal improvement in the absolute value of correlation with the evolving residual.
The soft thresholding parameter in the implicit formulation is set to $\alpha=0.1$ and $\alpha=0.2$ for  $\eta=0.1$ and $\eta=0.9$, respectively.

\begin{figure}[htb]
\centering
\includegraphics[width=\textwidth]{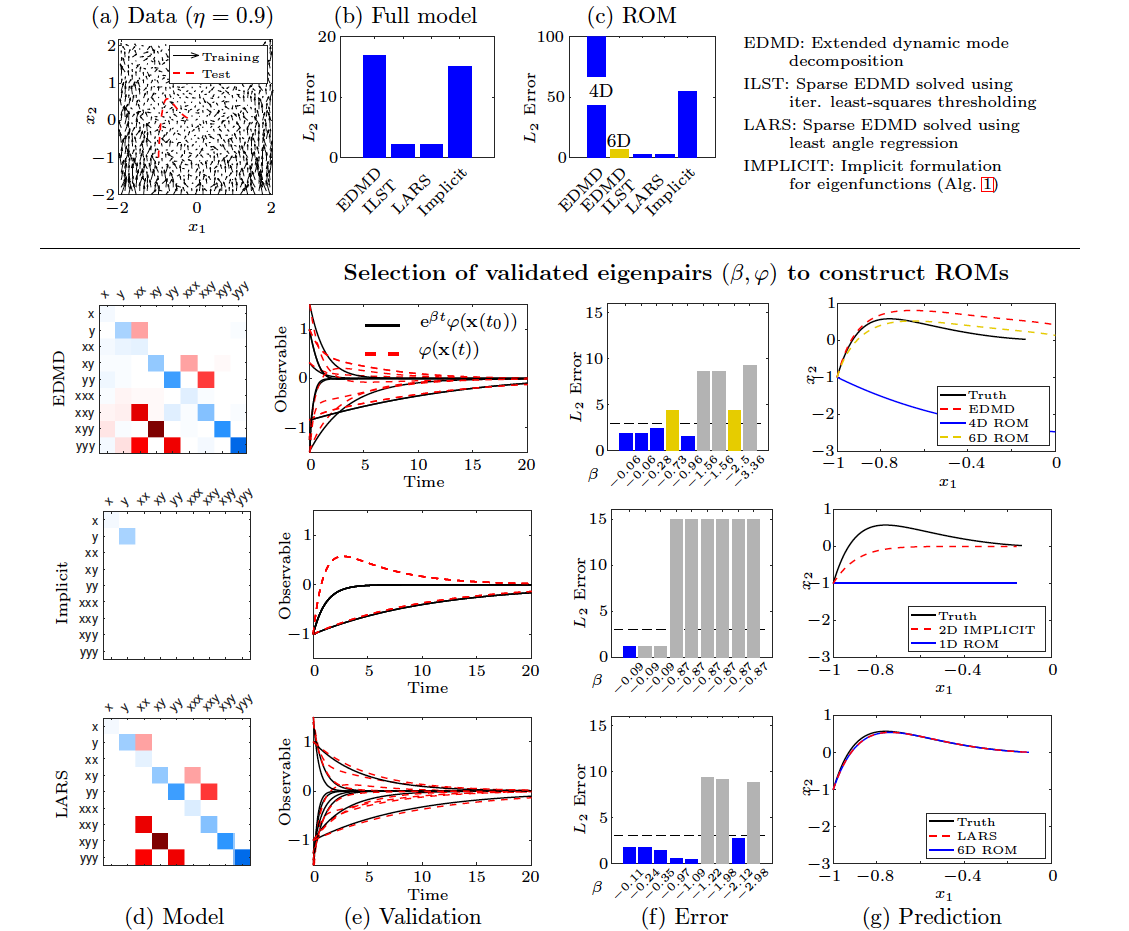}
	\caption{
	\RefereeOne{
		Model comparison for the slow manifold system as in Fig.~\ref{Fig:SLOWMANIFOLD:eta01} but with higher measurement noise level $\eta=0.9$. The ROM model error displayed in (c) is cut at $100$.
		For details see Fig.~\ref{Fig:SLOWMANIFOLD:eta01}.
	}
	}
	\label{Fig:SLOWMANIFOLD:eta09}
	\vspace{-0.in}
\end{figure}
The results are displayed in Figs.~\ref{Fig:SLOWMANIFOLD:eta01} and ~\ref{Fig:SLOWMANIFOLD:eta09} for the two noise cases.
We observe that with increasing noise level, the EDMD state transition matrix becomes denser due to overfitting (see (d)).
In contrast, the implicit approach and $L_1-$regularized EDMD (termed `LARS' or `ITLS' in the following) yield sparse matrices, which is more apparent for large noise magnitude.
Note that the state matrix for the implicit formulation is reconstructed using the identified eigenfunctions on the set of candidate functions for visualization purposes, 
i.e. there are many zero entries and candidate functions that do not contribute to the actual dynamics.
Lightly damped eigenfunctions can be recovered for all approaches. 
Eigenfunctions, or specifically eigenpairs $(\beta,\varphi(\bx))$, are deemed as recovered, if the prediction error falls below a threshold (marked as dashed line in (f)).
The prediction error is computed by evaluating the eigenfunction on a test trajectory and comparing it with the evolution as predicted by the associated eigenvalue; the evolution is displayed in (e).
These \emph{good} eigenpairs (colored as blue bars in (f)) are then used to construct reduced-order models (ROMs), that are by design linear and closed (subject to small errors). 
Note that solution from the implicit formulation may not be unique, i.e. there are several identical eigenpairs.
However, only unique eigenpairs are used to construct the reduce-order model. This is the reason why, e.g. in Fig.~\eqref{Fig:SLOWMANIFOLD:eta01}(f), the error for all eigenpairs falls below the threshold, but only unique ones (marked in blue) are selected to construct the ROM.

For small noise magnitude (Fig.~\ref{Fig:SLOWMANIFOLD:eta01}), the identified models and ROMs perform all well, except the EDMD-based ROM, as one crucial eigenfunction associated with $\beta = -0.19$ falls above the threshold and is not selected. The eigenvalue is correct, but not the associated eigenfunction.
If the eigenfunction is included in the model (yellow bar in (f) and yellow dashed line in (g)), i.e. the model is constructed from blue and yellow marked eigenpairs, it yields a similar performance as the full-state EDMD model.
It can also be observed that even in the low noise setting sparsification yields improvements. 
We note that the full-state model and ROM obtained from the implicit formulation are identical, as all discovered eigenpairs fall below the threshold (in (f)) and only unique pairs are selected to construct either model.
We note that the 4D-ROM from the implicit formulation can achieve better accuracy with one dimension lower compared with EDMD and the sparse EDMD model identified using LARS. 
Moreover, while all eigenpairs from the implicit formulation are accurate, both EDMD and LARS yield also non-physical eigenfunctions. 
The overall $L_2$ prediction error on the test trajectory in the original state is summarized for the full-state models and ROMs in (b) and (c), respectively.

For higher noise levels (Fig.~\ref{Fig:SLOWMANIFOLD:eta09}), the performance differences become more apparent.
LARS still yields several more accurate eigenfunctions than EDMD, so that LARS and the LARS-ROM significantly outperform EDMD and the EDMD-based ROM.
The performance of EDMD does not improve here by truncating the SVD when computing the pseudo-inverse.
The implicit formulation is more noise sensitive, as it is searching for eigenfunctions in the nullspace of a matrix~\ref{Eq:SparseKoopman}, and is only able to discover one accurate eigenfunction. 
Interestingly, despite this caveat, the 2D model is able to outperform EDMD and the EDMD-based ROM.
The implicit model is unable to predict the exact transient behavior; however, it accurately predicts the convergence to the steady state.
The implicit ROM model, however, is insufficient with just one eigenfunction $\varphi(\bx)=x_1$ and is unable to predict the evolution of the second state $x_2$.
When selecting one or two additional eigenfunctions for the EDMD-based ROM (in yellow  shown in (g) for two additional eigenfunctions, i.e. in total six eigenfunctions marked by blue and yellow in (f)), 
the prediction performance improves considerably, although it is still unable to predict the steady-state behavior.

Summarizing, EDMD is suffering from overfitting and sparsity-promoting formulations can yield significant performance enhancements.
Model validation in terms of eigenfunctions is a crucial step and can be used to construct better performing and lower-dimensional reduced-order models.
More generally, these results also demonstrate the importance of not just selecting \emph{good}, but also the \emph{right} eigenfunctions.

\subsection{Control design}
We now demonstrate control in intrinsic Koopman coordinates for the controlled system:
%
\begin{align}\label{Eqn:SlowManifold}
\frac{\mathrm d}{\mathrm dt}\begin{bmatrix} x_1 \\ x_2 \end{bmatrix} = \begin{bmatrix} \mu x_1 \\ \lambda (x_2 - x_1^2)\end{bmatrix} + \bB\, u
\end{align}
where the control vector is $\bB \in\mathbb{R}^2$.
%
This system can be represented as a finite-dimensional, linear system in a special choice of observable functions, making it amenable to optimal control~\cite{Brunton2016plosone}.  
KRONIC in intrinsic coordinates provides a powerful alternative if the system does not allow for a fully controllable, linear representation.

The system	exhibits slow and fast dynamics for $ \vert\lambda\vert \ll \vert\mu\vert$ and has a single fixed point at the origin.
This nonlinear system can be embedded in a higher-dimensional space 
$(y_1,y_2,y_3) = (x_1,x_2,x_1^2)$ where the unforced dynamics form a closed linear system in a Koopman-invariant subspace:
\begin{align}\label{Eqn:SlowManifold_MeasSys}
\frac{\mathrm d}{\mathrm dt}\begin{bmatrix} y_1 \\ y_2 \\ y_3\end{bmatrix} = \underbrace{\begin{bmatrix} \mu & 0 & 0\\ 0 & \lambda & -\lambda \\ 0 & 0 & 2\mu\end{bmatrix}}_{\bK}
\begin{bmatrix}
y_1 \\ y_2\\ y_3
\end{bmatrix} + 
\underbrace{
\begin{bmatrix}
1 & 0\\0 & 1\\2y_1 & 0
\end{bmatrix}
\bB}_{=\bB_y} \,u.
\end{align}
However, $\bB_y$ may be a function of $\by$, and hence of state $\bx$, depending on the specific choice of $\bB$.
Koopman eigenfunctions of the unforced system, i.e. $\bB \equiv [0\,\,0]^T$, 
are $\varphi_{\mu} = x_1$ and $\varphi_{\lambda} = x_2 - b x_1^2$ with $b = \frac{\lambda}{\lambda-2\mu}$  with eigenvalues $\lambda$ and $\mu$, respectively. 
These eigenfunctions remain invariant under the Koopman operator $\bK$ and can be interpreted as intrinsic coordinates.  
Note that $\varphi_{p\beta} := \varphi_{\beta}^p$ are also Koopman eigenfunctions with eigenvalue $p\beta$ for $p\in\mathbb{N}$ (and $p\in\mathbb{Z}$ for non-vanishing $\varphi_{p\beta}$).

The dynamics of the Koopman eigenfunctions are affected by the additional control term $\bB {\not =} [0\,\,0]^T$ according to (see also Eq.~\eqref{Eqn:KRONIC:MIMO})
\begin{align}\label{Eq:KoopmanEigenfunctionSys}
	\frac{\mathrm d}{\mathrm dt} \boldsymbol{\varphi}
= \begin{bmatrix}
	\mu & 0 & 0\\
	0     & \lambda & 0 \\
	0 & 0 & 2\mu
	\end{bmatrix}\boldsymbol{\varphi}
	+ \begin{bmatrix}
	1 & 0\\
	-2 b x_1 & 1\\
	2 x_1 & 0
	\end{bmatrix} \cdot  \bB\, u.
\end{align}
Here, the first term represents the unforced, uncoupled, linear dynamics of the eigenfunctions and the second term a possibly state-dependent control term $\nabla \boldsymbol{\varphi}\cdot\bB$, that incorporates the effect of control on each eigenfunction. 

The controller shall stabilize the unstable fixed point at the origin if either $\mu$ or $\lambda$ are unstable. 
The control objective is to minimize the quadratic cost function
\begin{equation}\label{SlowManifold_Cost_x}
J_x =\int_{0}^{\infty}\, \bx^T \bQ \bx + R\,u^2\,\mathrm dt
\end{equation}
with $\bQ = \left[\begin{smallmatrix}1 & 0\\ 0 & 1\end{smallmatrix}\right]$ and $R = 1$, weighing state and control expenditures equally. 
Analogously, 
we can define a cost function in observable functions,
\begin{eqnarray}\label{Eqn:SlowManifold:Cost_Jy}
J_y = \int_{0}^{\infty}\, \by^T \bQ_y \by + R\,u^2\,\mathrm dt, \quad \bQ_y = \left[\begin{smallmatrix}
\bQ  & & 0\\ & & 0 \\ 0 & 0 &0
\end{smallmatrix}\right],
\end{eqnarray}
and in intrinsic coordinates,
\begin{eqnarray}\label{Eqn:SlowManifold:Cost_Jphi}
J_{\varphi} = \int_{0}^{\infty}\, \bvarphi^T \bQ_{\varphi} \bvarphi + R\, u^2\,\mathrm dt, \quad \bQ_{\varphi} = \left[\begin{smallmatrix}
\bQ & & 0\\ & & b \\ 0 & b &b^2
\end{smallmatrix}\right].
\end{eqnarray}
Here, $\bQ_{y}$ and $\bQ_{\varphi}$ are chosen to yield the same cost in $\bx$. 
Linear optimal control is then directly applied to~\eqref{Eq:KoopmanEigenfunctionSys} to derive the control law, which is then incorporated in~\eqref{Eqn:SlowManifold}.
The controller is linear in $\by$ and $\boldsymbol\varphi$ and yields a nonlinear controller in the state $\bx$:
	\begin{align*}
	u_y =& - \bC_y \by  
	= - \begin{bmatrix} C_{y,1} & C_{y,2}\end{bmatrix} \begin{bmatrix}
	x_1\\x_2 \end{bmatrix} - C_{y,3}\, x_1^2\\
	u_{\varphi} =& - \bC_{\varphi} \boldsymbol{\varphi} = -\begin{bmatrix} C_{\varphi,1} & C_{\varphi,2}\end{bmatrix}\begin{bmatrix} x_1 \\ x_2\end{bmatrix} - (C_{\varphi,3}-bC_{\varphi,1})x_1^2, 
	\end{align*}	
where $\bC_y\in\mathbb{R}^{1\times 3}$ and $\bC_{\varphi}\in\mathbb{R}^{1\times 3}$
are the control gain vectors in observable or intrinsic coordinates, respectively.

\begin{figure}[tb]
	\centering
	\includegraphics[width=\textwidth]{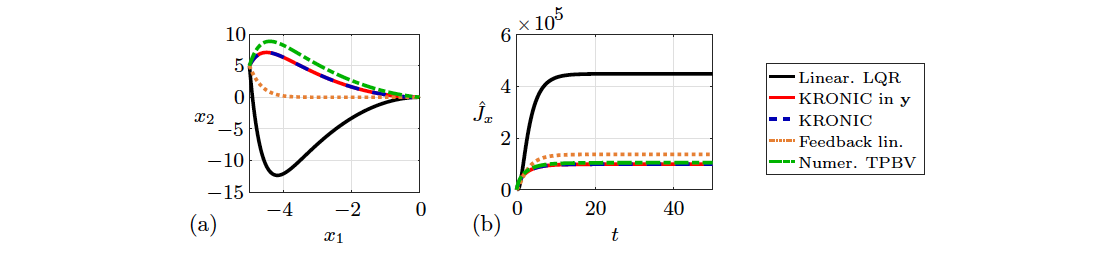}
	\caption{LQR for $\bB = [0\;\; 1]^T$, $\mu = -0.1$ and $\lambda=1$ using standard linearization, truncated Koopman in $\by$ and $\boldsymbol\varphi$ (KRONIC) compared with the solution of the nonlinear control problem (TPBV) and feedback linearization:
		(a) phase plot and (b) cost evolution. 	
		KRONIC is outperforming all other approaches.}
	\vspace{-.08in}
	\label{Fig:SlowManifold_B01}
\end{figure}
\subsection{Stabilization with unstable fast dynamics}
First, we consider the system  with $\mu=-0.1$ and $\lambda=1$ with an unstable $x_2$ direction. 
The control vector is $\bB=[0\;1]^T$, 
resulting in a constant vector in $\by$ or $\boldsymbol\varphi$ coordinates, $\bB_y=\bB_{\varphi} = [ 0 \; 1 \; 0 ]^T$.  
Note that the first direction is uncontrollable, but also stable.

Standard LQR results are compared (see Fig.~\ref{Fig:SlowManifold_B01}) for the linearized dynamics, truncated Koopman system in $\by$, truncated Koopman system in $\boldsymbol\varphi$ (KRONIC), as well as with feedback linearization~\cite{khalil1996book} ($u_{FL}=\lambda x_1^2 - \bC_{FL}\bx$) and with numerically solving the nonlinear control problem as a two-point boundary value problem (TPBV), with performance evaluated in terms of the cumulative cost $\hat{J}_x^t = \sum_{\tau=0}^t\,J_x(\tau)$. 
Both controllers, in observable functions and intrinsic coordinates,  
achieve the same performance and outperform linearized LQR and feedback linearization. 
The results for the truncated Koopman system in observables correspond to those presented in~\cite{Brunton2016plosone}. 
There is no difference between those results and the control results of the system in intrinsic coordinates, as these systems are connected via an invertible linear transformation.
One advantage of a formulation in intrinsic coordinates will become apparent in the next case, where the stable and unstable directions are reversed.  

\begin{figure}[tb]
\centering
\includegraphics[width=\textwidth]{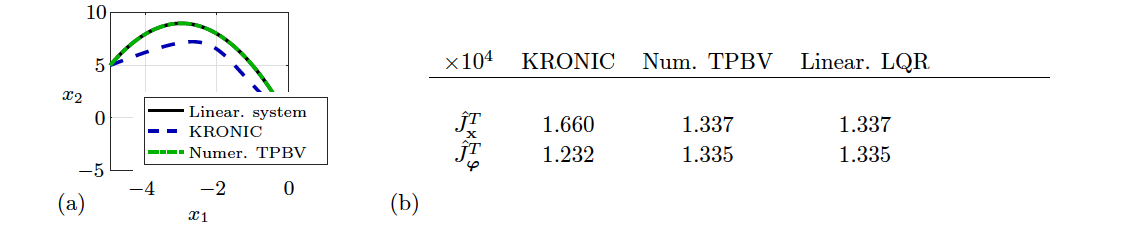}
	\caption{
		Control for $\bB = [1\;0]^T$, $\mu = 0.1$ and $\lambda=-1$ comparing KRONIC~\eqref{Eqn:SlowManifold:KoopmanIntrinsicCoordinates_Case2} by solving the SDRE, 
					Num. TPBV, and LQR on the linearized dynamics: 
					(a) phase plot, and (b) cumulative control performance over $T=50$ time units.
					Note that the cost for KRONIC can not be directly compared with TPBV and linear. LQR, as the controller in the latter two is optimized with respect to $J_x$, while the former is optimized for $J_{\boldsymbol{\varphi}}$.}
	\label{Fig:SlowManifold_B10}
\end{figure}

\subsection{Stabilization with unstable slow dynamics}
We now consider~\eqref{Eqn:SlowManifold} where the stable and unstable directions are reversed, i.e.\ $\mu=0.1$ and $\lambda=-1$. The control input now affects the first state $x_1$ with $\bB = [1\;0]^T$, otherwise this state is uncontrollable.
As elaborated in~\cite{Brunton2016plosone}, in this case the linear system in observables~\eqref{Eqn:SlowManifold_MeasSys} will become nonlinear in the control term and, more importantly, will become unstabilizeable as the third state $y_3$ has a positive eigenvalue $2\mu$. 
Analogously, the Koopman system in intrinsic coordinates has an uncontrollable, unstable direction in $\varphi_{2\mu}$.
However, the dynamics of the Koopman eigenfunctions are uncoupled, thus the third direction can be discarded and the controller is developed in the controllable subspace:
\begin{equation}\label{Eqn:SlowManifold:KoopmanIntrinsicCoordinates_Case2}
\frac{\mathrm d}{\mathrm dt} 
\begin{bmatrix}
\varphi_{\mu}\\
\varphi_{\lambda}
\end{bmatrix}
= 
\begin{bmatrix}
\mu & 0 \\
0     & \lambda  
\end{bmatrix} 
\begin{bmatrix}
\varphi_{\mu}\\
\varphi_{\lambda}
\end{bmatrix}
+ 
\begin{bmatrix} 1 \\ -2 b x_1
\end{bmatrix} u.
\end{equation}
%
Note that the third direction $\phi_{2\mu}(x)$ is a harmonic of $\phi_{\mu}(x)$, i.e., $\phi_{2\mu}=\phi_{\mu}^2$.  
Thus, these two directions may not be independently controllable with a single input.
Here, the controller for $\bx$ is determined by solving the SDRE (see Sec.~\ref{Sec:KRONIC:OptimalControl}) at each time instant to account for the nonlinear control term.
In a truncated Koopman eigenfunction system like in this example, the weights in $J_{\boldsymbol{\varphi}}$ can generally not be modified to directly replicate the cost $J_{\bx}$.
Here, the weights for $J_{\bx}$ are $\bQ = \texttt{eye}(2)$ and $R=1$, and for $J_{\boldsymbol{\varphi}}:=J(\boldsymbol{\varphi},\bu)$	are $\bQ_{\boldsymbol{\varphi}} = \bQ$ and $R_{\boldsymbol{\varphi}}=4R$.
The choice for $R_{\boldsymbol{\varphi}}$ ensures a fair comparison by enforcing the same applied energy input for all methods.
	
Performance results are summarized in Fig.~\ref{Fig:SlowManifold_B10}.
It is also possible to combine KRONIC with MPC allowing for more general objective functions. 
Then, the control could be formulated in terms of $J_{\bx}$ by computing the inverse $\boldsymbol{\varphi}^{-1}:\mathbb{R}^r\rightarrow \mathbb{R}^n$ if it exists or estimating $\bx$ from $\boldsymbol{\varphi}$ using, e.g., multidimensional scaling as in~\cite{kawahara2016nips}.
Control in Koopman intrinsic coordinates
allows one to discard uncontrollable, unstable directions, for which standard control toolboxes such as Matlab's \texttt{lqr} fail. 
Note that feedback linearization fails in this case too: The control law is of the form $u_{\mathrm{FL}} = x_1^{-1} \bC_{\mathrm{FL}} \bx$. As the system approaches the origin, the control input becomes unboundedly large. 

%% file: HamiltonianSystems.tex
\begin{figure}[tb]
	\centering
	\begin{overpic}[width=1\textwidth]{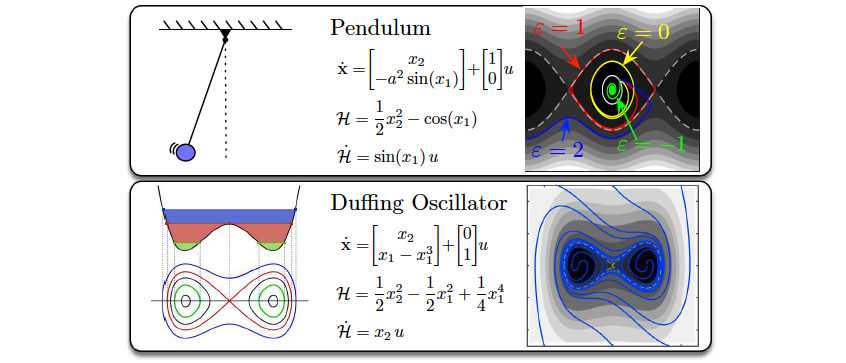}
	\end{overpic}
	\vspace{-.1in}
	\caption{KRONIC demonstrated for several Hamiltonian systems.}
	\label{Fig:HamiltonianSystems:Overview}
	\vspace{-.1in}
\end{figure}

Conserved quantities, such as the Hamiltonian, are Koopman eigenfunctions associated with the eigenvalue $\lambda=0$. 
Hamiltonian systems represent a special class of systems for which we can easily discover a Koopman eigenfunction directly from data.

The dynamics of a Hamiltonian system are governed by
\begin{equation}\label{Eqn:HamiltonianSystem:EvolutionEquation}
	\frac{\mathrm d}{\mathrm dt}\bq =\frac{\partial \mathcal{H}}{\partial \bp},\quad
	\frac{\mathrm d}{\mathrm dt}\bp = -\frac{\partial \mathcal{H}}{\partial \bq},
\end{equation}
where $\bq$ and $\bp$ are the generalized state and momenta vectors, respectively.  
The Hamiltonian $\mathcal{H} = \mathcal{H}(\mathcal{\bq,\bp})$ considered here is time-independent, representing the conserved energy in the system. 
Trajectories of the system evolve on constant energy hypersurfaces $\{(\bq, \bp):\, \mathcal{H}(\bq, \bp) = \mathcal{E}\}$, 
which may be interpreted as oscillatory modes.
Thus, energy level stabilization is a form of oscillation control and corresponds to stabilizing invariant manifolds in phase space.
Also nonlinear fixed point stabilization may correspond to stabilizing a particular value of the Hamiltonian energy.

Consider the nonlinear, control-affine Hamiltonian system
$
\frac{\mathrm d}{\mathrm dt}\bx = {\boldsymbol{\it f}}(\bx) + \bB \bu
$ 
where ${\boldsymbol{\it f}} = [ \partial\mathcal{H}/\partial\bp\;-\partial\mathcal{H}/\partial\bq]^T$,
with state vector $\bx = [\bq \,\,\bp]^T\in\mathbb{R}^n$, multi-channel input vector $\bu\in \mathbb{R}^{q}$, and constant control matrix $\bB \in\mathbb{R}^{n\times q}$. 
We may develop the control directly for the eigenfunction equation
\begin{equation}
\frac{d}{dt}\varphi = 0\cdot\varphi + {\nabla_{\bx}\varphi\cdot\bB}\bu := \bB_{\varphi}\bu,
\end{equation}
where $\varphi = \mathcal{H}$.
This equation also represents the energy conservation law of the system: A change in the energy, i.e.\ the Hamiltonian,  corresponds to the external supplied work via $\bu$.
\RefereeTwo{
The infinite-time horizon cost function to be minimized is 
$
J =  \lim_{t\rightarrow\infty}\frac{1}{2} \int_0^{t} Q\left(\mathcal{H}(\bx(t))\right)^2 + 
 \bu^T(t) \bR \bu(t) \,dt	$
with scalar $Q$ penalizing energy deviations and $\bR$ penalizing the cost expenditure.
}
Assuming a single input $u$, the control law is given by
$
u = - \mathrm{sign}(B_{\mathcal{H}}) \sqrt{{Q}/{R}}\mathcal{H}(\bx)
$ 
feeding back the current energy level $\mathcal{H}(\bx)$. The ratio $Q/R$ determines how aggressive the controller is. 
A more aggressive controller, with $Q>R$, leads to a faster but also more costly convergence to the desired state, and vice versa.
For the specific case  $Q=R$, 
the control law reduces further to $u = - \mathrm{sign}(B_{\mathcal{H}}) \mathcal{H}(\bx)$.
Note that this feedback control law is linear in the Hamiltonian function, but nonlinear in the state $\bx$.
In the following, we demonstrate the control approach for several Hamiltonian systems by solving the SDRE; an overview is provided in Fig.~\ref{Fig:HamiltonianSystems:Overview}, where colored curves represent Koopman controlled trajectories.
We assume $Q=R=1$ for all examples.

\subsection{Frictionless pendulum}
A frictionless pendulum given by a negative cosine potential $V(x) = -\cos(x_1)$  is steered 
\RefereeTwo{
towards different energy levels $\mathcal{H}^{REF}(\bx) = \mathcal{E}\,\forall\bx$, $\mathcal{E}\in\{-1,0,1,2\}$, with feedback law  $u= - C(\mathcal{H}(\bx)-\mathcal{H}^{REF}) = - C(\mathcal{H}(\bx)-\mathcal{E})$.
}
Minimizing the Hamiltonian drives the system to the equilibrium point at the origin with $\mathcal{E}=-1$ (green curve in Fig.~\ref{Fig:HamiltonianSystems:Overview}).  
Koopman control has improved performance over an LQR controller based on linearized dynamics near the center.  

\begin{figure}[tb]
\centering
\includegraphics[width=\textwidth]{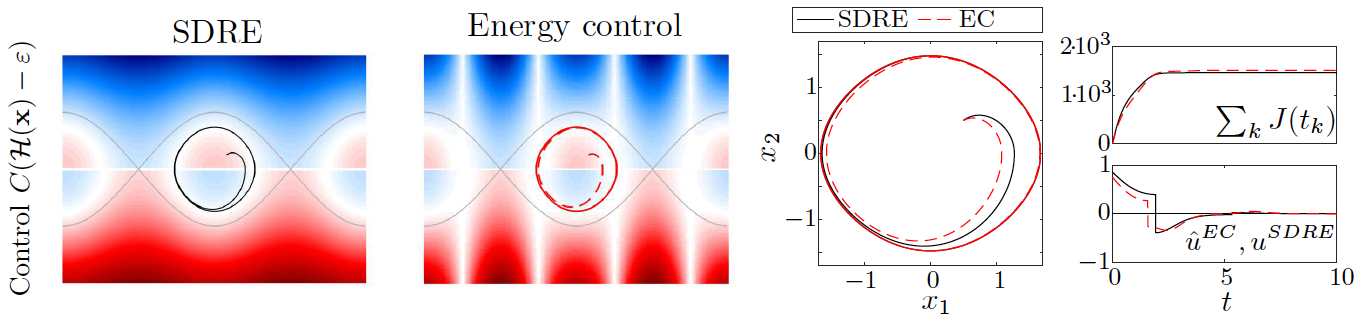}
	\caption{Comparison of the energy control control strategy~\cite{Astrom2000autom} and SDRE for the pendulum with $\varepsilon=0.1$. The color coding in the phase plots depicts the respective state-dependent control}
	\label{Fig:HamiltonianSystems:Pendulum}
	\vspace{-.1in}
\end{figure}

We further compare the controller resulting from solving the state-dependent Riccati equation with an energy-based control (EC)~\cite{Astrom2000autom}, which is developed from physical considerations:
\begin{subequations}
\begin{align}
SDRE:& \; \bB = \begin{bmatrix}0\\1\end{bmatrix},
\quad \dot{\mathcal{H}} = x_2 u,\quad u =  - \mathrm{sign} (x_2)(\mathcal{H}(\bx)-\mathcal{E})\\
EC:& \; \bB = \begin{bmatrix}0\\-\cos(x_1)\end{bmatrix},
\; \dot{\mathcal{H}} = -x_2\cos(x_1) u, \; u = - \mathrm{sign} (-x_2\cos(x_1))(\mathcal{H}(\bx)-\mathcal{E})\label{Eqn:Pendulum:EC}
\end{align}		
\end{subequations}
If the functional representation of the Koopman eigenfunction, or specifically here the Hamiltonian energy function, is identified from data as in KRONIC, the feedback gain can be precomputed and does not need to be computed on-line.
In Fig.~\ref{Fig:HamiltonianSystems:Pendulum}, the
state-dependent feedback law is shown for both controllers as color-coded phase plots, where negative values are depicted by blue and positive values by red. The extrema are identical.
The control input from the EC controller is here slightly modified compared with~\eqref{Eqn:Pendulum:EC} by moving the term $-\cos(x_1)$ into $u$ itself. Then the effective control and performance are displayed in terms of $\hat{u} = -\cos(x_1) u$ for a fair comparison; otherwise, the EC controller would appear less performing.
The pendulum dynamics are locally controllable as long as $u$ does not vanish, which is achieved
if $x_1{\not =}\pi/2$ and $x_2 {\not =}0$ for the EC controller and if $x_2 {\not =}0$ for the SDRE controller.
Both controllers yield a bang-bang strategy and successfully steer the system to the desired energy level, though the SDRE controller is performing slightly better.

\subsection{Duffing system}
Certain damped or forced oscillators are described by the Duffing equation (Fig.~\ref{Fig:HamiltonianSystems:Overview}).
The system is steered towards the energy level $\mathcal{E}=0$, which corresponds to the separatrix cycle (yellow dashed lines) yielding a periodic solution.
The origin is a saddle point leading to a homoclinic orbit defined by $x_2^{*} =  \pm x_1^{*}\sqrt{1-{1} (x_1^{*})^{2}/2}$ for $x_1^{*} \leq \pm\sqrt{2}$.
It is possible to stabilize the center fixed points by commanding a lower reference energy; however, because of symmetry in the system, both fixed points are indistinguishable in eigenfunction coordinates.  
This illustrates a fundamental \emph{uncertainty} associated with Koopman eigenfunction control.

\begin{figure}[tb]
\centering
\includegraphics[width=\textwidth]{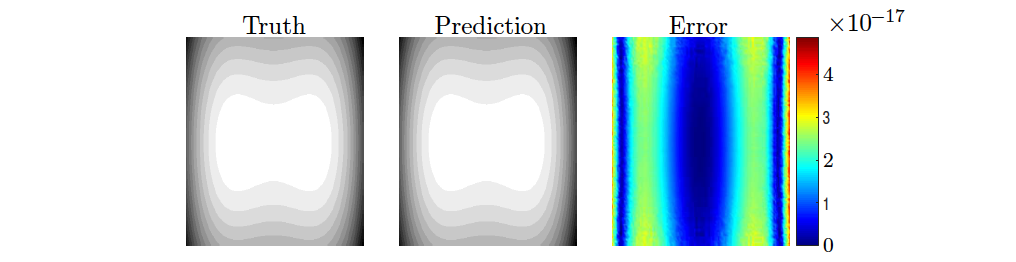}
	\caption{Prediction over the state space and error using the discovered Hamiltonian for the Duffing system trained on a single trajectory.}
	\label{Fig:Duffing:PredictionError}
\vspace{-0.1in}
\end{figure}

For the eigenvalue $\lambda=0$,~\eqref{Eq:SparseKoopman} becomes $-\bGamma(\bX,\dot{\bX})\bxi = {\bf 0}$, and hence a sparse $\bxi$ is sought in the null-space of $-\bGamma(\bX,\dot{\bX})$.
Polynomials up to fourth order are employed to construct a library of candidate functions in~\eqref{Eq:Theta}-\eqref{Eq:Gamma}.
A single time series for $t\in[0,10]$ with time step $\Delta t = 0.001$ starting at the initial state $\bx_0 = [0,\; -2.8]^T$ is collected.
Thus, each row in~\eqref{Eq:Theta} and~\eqref{Eq:Gamma} corresponds to a time instant of the trajectory.
The prediction and associated error of the identified eigenfunction are displayed in Fig.~\ref{Fig:Duffing:PredictionError}.
Here, the energy is predicted over the full state space using the identified eigenfunction, the Hamiltonian energy function, which is trained from a single trajectory.
The magnitude of the error is very small of $\mathcal{O}(10^{-17})$. However, the eigenfunction evaluate on a trajectory would oscillate with a tiny amplitude around the true energy level. 

\begin{figure}[tb]
\centering
\includegraphics[width=\textwidth]{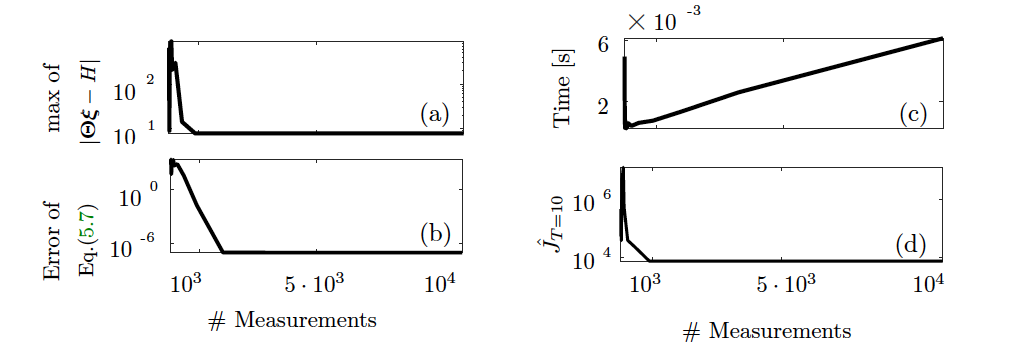}
	\caption{Dependency in increasing number of measurements for the Duffing system: 
			(a) error in the Koopman eigenfunction,
			(b) error in~\eqref{Eq:SparseKoopman} with estimated $\boldsymbol{\xi}$, 
			(c) computational time to solve for $\boldsymbol{\xi}$ in~\eqref{Eq:SparseKoopman}, and
			(d) cumulative control performance $\hat{J}$ over $T=10$ time units.}
	\vspace{-.05in}
	\label{Fig:DuffingStats}
\end{figure}
The error of the regression problem, the computational time for identifying the eigenfunction, and the control performance (steering towards $\mathcal{E}=0$) for an increasing number of measurements in the estimation step are displayed in Fig.~\ref{Fig:DuffingStats}.
The identified Koopman eigenfunction with $\lambda=0$ from 1792 measurements (kink in Fig.~\ref{Fig:DuffingStats}(b)) is
$
\varphi(\bx)= \left[
\begin{smallmatrix}
x_1^2 \; x_2^2 \; x_1^4
\end{smallmatrix}\right]
\left[\begin{smallmatrix}
-\frac{2}{3} \; \frac{2}{3} \; \frac{1}{3}
\end{smallmatrix}\right]^T$ with error $\mathcal{O}(10^{-8})$. 
This eigenfunction represents a perfect recovery of the Hamiltonian up to a scaling, as a Hamiltonian multiplied by a constant scalar is also a conserved quantity.
Using a larger time step of $\Delta t=0.05$, 56 measurements are sufficient to learn the eigenfunction with error $\mathcal{O}(10^{-6})$.

%% file: AsymmetricPotentialWell.tex
\begin{figure}[tb]
\centering
\includegraphics[width=\textwidth]{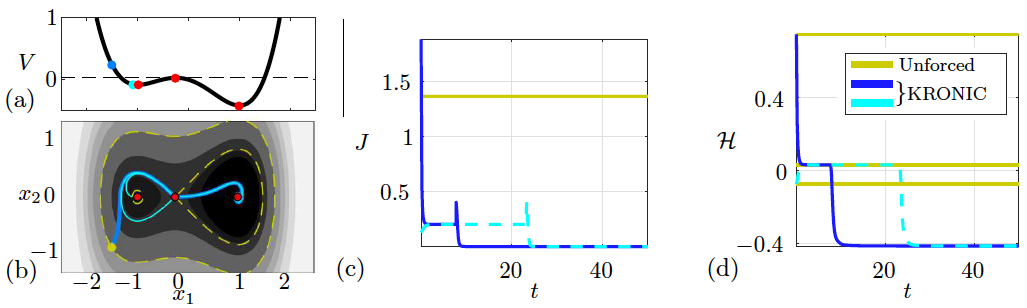}
	\caption{
		Switching control strategy based on KRONIC and the homoclinic orbit to jump between wells:
		(a) potential function $V(x_1)$ showing initial conditions (blue and cyan) and extrema (red), 
		(b) phase plot with unforced (yellow) and controlled trajectories (blue and cyan), 
		(c) instantaneous cost function and (d) total energy.}
		\vspace{-.13in}
	\label{Fig:AsymPotentialWell:Results}
\end{figure}

A Koopman eigenfunction represents a topography over the state; e.g., the Hamiltonian function depicts the energy landscape of the system.
Trajectory control of a set of particles based on a single Koopman eigenfunction is driven by the difference between the current and desired value in this topography.
While the Koopman eigenfunction is a global, linear representation of the system, the control of the particle is local, e.g. by solving a suboptimal SDRE.
This is illustrated for a particle in an asymmetric double potential
$V(x_1) = \frac{1}{4}x_1^4 - \frac{1}{2}x_1^2 - \frac{a}{3} x_1^3 + a x_1$ 
with $a=-0.25$ (Fig.~\ref{Fig:AsymPotentialWell:Results}(a)). 
The Hamiltonian is $\mathcal{H} = x_2^2/2 + V(x_1)$ and the dynamics are
\begin{equation}
\frac{\mathrm d}{\mathrm dt} \begin{bmatrix} x_1\\ x_2 \end{bmatrix} 
= 
\begin{bmatrix}
x_2 \\ -x_1^3 + ax_1^2 + x_1 - a
\end{bmatrix}
+ 
\begin{bmatrix} 1& 0 \\ 0 & 1 \end{bmatrix} \bu.
\end{equation}
For an initial condition in the left well (blue dots in Fig.~\ref{Fig:AsymPotentialWell:Results}(a)) 
the controller will fail to steer the state to the fixed point $\bx^{*} = [1\,\,0]^T$ in the center of the right well 
as the trajectory will become trapped in the bottom of the left well.  
Instead, the controller must first increase the energy level to the saddle transition, and after the trajectory passes to the right basin, the energy can be decreased further.  

We propose a switching control strategy that exploits the Koopman eigenfunctions to transport particles between basins of different heights associated with different fixed points.
In particular, the following control strategy steers particles from the left to the right basin:
\begin{equation*}
\bu = \begin{cases}
-\bC(\bx) (\mathcal{H}(\bx)-\mathcal{H}([a,0])) & \text{if } \mathcal{H}(\bx) < \mathcal{H}([a,1]),\\
0 & \text{if } \mathcal{H}(\bx) = \mathcal{H}([a,1]) \text{and }  x_1\leq a,\\
-\bC(\bx) (\mathcal{H}(\bx)-\mathcal{H}([1,0])) & \text{if } x_1> a.
\end{cases}
\end{equation*}
A particle on an energy level lower than $\mathcal{H}([a,1])$, associated with the saddle point, is first steered onto a trajectory with slightly higher energy than the homoclinic orbit connecting the two basins.
On this orbit, control is turned off and the particle travels to the right well exploiting the intrinsic system dynamics. 
As soon as it passes the saddle point, control is turned on again directing it to the lowest energy level $\mathcal{H}([1,0])$ at the desired fixed point. 	
The controller is demonstrated for two initial conditions, as shown in Fig.~\ref{Fig:AsymPotentialWell:Results}(b-d), driving both to the desired energy level. 

This controller can be fully derived from data: First, relevant Koopman eigenfunctions can be identified from data, as shown in Sec.~\ref{Sec:Regression}. 
By analyzing roots and extrema of the eigenfunction corresponding to $\lambda=0$, equilibrium and saddle points can be identified.
The homoclinic and heteroclinic orbits associated with the saddles can be used as natural transit paths between basins. 
In each basin, eigenfunction control drives the system to the desired state.
Future applications for this control strategy include space mission design and multi-stable systems such as proteins.

%% file: DoubleGyre.tex
We now consider a high-dimensional, spatially evolving, non-autonomous dynamical system, with time-dependent Koopman eigenfunctions.
The periodically driven double gyre flow models the transport between convection cells in the Rayleigh-B\'enard flow due to lateral oscillations, yielding a simple model for the gulf stream ocean front~\cite{Solomon1988pra}.
We employ here the same parameters as in Shadden et al.'s seminal work on Lagrangian coherent structures~\cite{Shadden2005physica}.

The time-dependent stream function is
\begin{subequations}
	\begin{align}
		\Psi(x,y,t) &= A\,\sin(\pi f(x,t))\sin(\pi y)\\
		\text{with}\quad f(x,t) &= \varepsilon \sin(\omega t)x^2 + (1-2\varepsilon\sin(\omega t))x 
	\end{align}
\end{subequations}
with $A=0.25$, $\omega = 2\pi$, and $\varepsilon=0.25$ on a periodic domain $[0,2]\times[0,1]$.
The velocity field ${\bf v} = [ v_x \,\, v_y]^T$ is given by
$v_x = -\frac{\partial\Psi}{\partial y}$ and
$v_y = \frac{\partial\Psi}{\partial x}$.
 
The control objective is to steer an ensemble of trajectories to a level set of the stream function.
This can be interpreted as the control of an ensemble of active drifters or autonomous gliders in the ocean, 
which drift due to hydrodynamic forces associated with $v_x$ and $v_y$. 
The dynamics of the $i$th drifter are
$\frac{\mathrm d}{\mathrm dt}\left[\begin{smallmatrix} x_i\\ y_i \end{smallmatrix}\right]
= \left[\begin{smallmatrix} v_x +\gamma_i \sin(\theta_i)\\ v_y+\gamma_i\cos(\theta_i)\end{smallmatrix}\right] 
= \left[\begin{smallmatrix} v_x \\ v_y\end{smallmatrix}\right] + \left[\begin{smallmatrix}1 & 0 \\ 0 & 1 \end{smallmatrix}\right] \bu$. 

For the autonomous and unforced flow with $\varepsilon=0$, the stream function is a Koopman eigenfunction associated with the eigenvalue $\lambda=0$. The forced system becomes:
\begin{equation}\label{Eqn:DoubleGyre:StreamFunDynamics_Autonomous}
\frac{\mathrm d}{\mathrm dt}{\Psi}(x,y) 
=
\underbrace{
	\begin{bmatrix}
	\frac{\partial\Psi}{\partial x} &  \frac{\partial\Psi}{\partial y}
	\end{bmatrix}
	\cdot
	\begin{bmatrix}
	-\frac{\partial\Psi}{\partial y} \\ \frac{\partial\Psi}{\partial x}
	\end{bmatrix}}_{=0}
+
\underbrace{
	\begin{bmatrix}
	\frac{\partial\Psi}{\partial x} &  \frac{\partial\Psi}{\partial y}
	\end{bmatrix}
	\cdot \bB}_{\bB_{\Psi}(x,y)}
\,\bu .
\end{equation}
Without control, the stream function $\Psi$ is conserved, as it is the negative of the Hamiltonian.  
The particles follow streamlines, which are isolines of the stream function. 

In the non-integrable case, with $\varepsilon>0$, the total derivative of the stream function is given by
\begin{equation}\label{Eqn:DoubleGyre:StreamFunDynamics}
\frac{\mathrm d}{\mathrm dt}{\Psi}(x,y,t) 
=  
A_{\Psi}(x,y,t) \Psi
+
\bB_{\Psi}(x,y,t)
\,\bu
\end{equation}
where the vanishing term in~\eqref{Eqn:DoubleGyre:StreamFunDynamics_Autonomous} is not displayed.
The first term in~\eqref{Eqn:DoubleGyre:StreamFunDynamics} arises from the time derivative $\frac{\partial}{\partial t}\Psi(x,y,t)$ and is reformulated into a linear-like structure in $\Psi$: \
$\partial \Psi/\partial t = A\pi \cos(\pi f(x,t))\sin(\pi y) (\partial f/\partial t) = \pi\tan^{-1}(\pi f(x,t))(\partial f/\partial t)\Psi = A_{\Psi} \Psi$. 
The second term is the time-dependent analogue of the corresponding term in~\eqref{Eqn:DoubleGyre:StreamFunDynamics_Autonomous}.

\begin{figure}[tb]
\centering
\includegraphics[width=\textwidth]{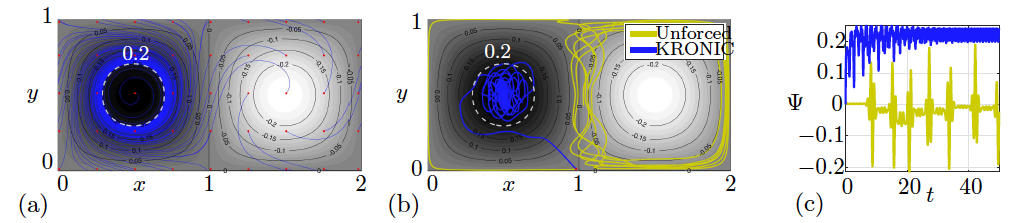}
	\caption{
		Controlled double gyre flow: (a) autonomous with $\varepsilon=0$ steering an ensemble of drifters (initial condition depicted by red dots) to the level $\Psi = 0.2$, 
		(b) a single drifter trajectory (unforced and controlled) in the non-autonomous double gyre flow with $\varepsilon=0.25$, and (c) corresponding stream function values.}
	\label{Fig:DoubleGyre}
\end{figure}

For both cases, $\varepsilon = 0$ and $\varepsilon>0$, a controller is developed for the stream function. The control is then applied to an ensemble of drifters to steer them towards the level set $\Psi = 0.2$.
As in the previous examples a quadratic cost function with $Q = 1$ and $\bR = (\begin{smallmatrix}
1 & 0 \\ 0 & 1 \end{smallmatrix})$ is considered.
In both cases, $\bB_{\Psi}$ and $A_{\Psi}$ depend on the state, and for $\varepsilon>0$ also on time. 
Thus, the state-dependent Riccati equation is solved at each point in space and time.
The controller successfully drives an ensemble of drifters distributed over the domain to the desired level, as shown in Fig.~\ref{Fig:DoubleGyre}(a). 
Trajectories are integrated using a $4$th-order Runge-Kutta scheme from $t\in[0,10]$.
Example trajectories of the non-autonomous system, with and without control, are presented in Fig.~\ref{Fig:DoubleGyre}(b-c).

Note that in the non-autonomous case, the reference isocurve $\Psi_{REF} = 0.2$ (white dashed in Fig.~\ref{Fig:DoubleGyre}(b)) oscillates from left to right while being periodically compressed and expanded in $x$-direction. The particles follow the moving isocurve resulting in a small oscillation around the desired value (see Fig.~\ref{Fig:DoubleGyre}(c)). 

Koopman eigenfunction control can be interpreted in two ways: 
(1) applying local control to internally driven swimmers or particles in an external field such as a fluid flow or magnetic field; or 
(2) driving the external field in which the swimmers or particles drift.
In the latter case, control drives the amplitude of the stream function at each point to the desired value.
For the double gyre flow with a constrained spatial domain, the spatially distributed gain may be precomputed.

%% file: Conclusions.tex
In summary, 
we extend the Koopman operator formalism to include actuation, and demonstrate how a nonlinear control problem may be converted into a bilinear control problem
in eigenfunction coordinates.  
Next, we have presented a data-driven framework to identify leading eigenfunctions of the Koopman operator and sparsity-promoting extensions to EDMD.
We find that lightly damped or undamped eigenfunctions may be accurately approximated from data via regression, as these eigenfunctions correspond to persistent phenomena, such as conserved quantities.  
Moreover, these are often the structures that we seek to control, since they affect long-time behavior.  
We have demonstrated the efficacy of this new data-driven control architecture on a number of nonlinear systems, including Hamiltonian systems and a challenging high-dimensional ocean mixing model.  
These results suggest that identifying and controlling Koopman eigenfunctions may enable significant progress towards the ultimate goal of a universal data-driven nonlinear control strategy.

\begin{figure}[tb]
\centering
\includegraphics[width=\textwidth]{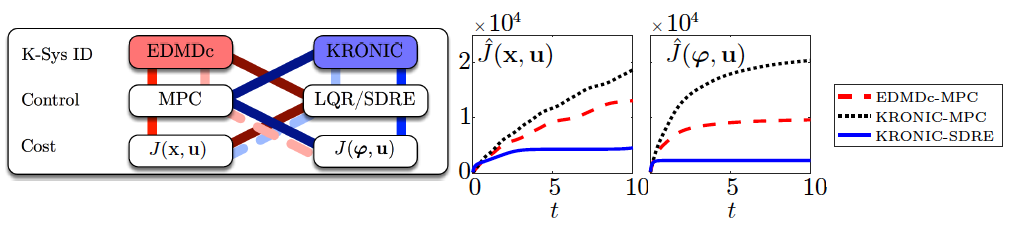}
	\caption{Comparison of Koopman system identification methods for control demonstrated for the Duffing system. Here, EDMDc-MPC requires knowledge of the Koopman eigenfunction. 
	The employed weights are $Q_H = 1$ and $R=1$ for $J_{\boldsymbol{\varphi}}$, and 
	$\bQ =\texttt{eye}(2)$ and $R=1$ for $J_{\bx}$ (instantaneous cumulative costs are shown here). The prediction and control horizon for MPC is in both cases $N=5$.
	}
	\vspace{-.12in}
	\label{Fig:KRONICvsEDMDc}
	\end{figure} 

\RefereeOne{
Finite-dimensional approximations of the Koopman operator are typically obtained as its projection onto a specified basis, which may suffer from a well-known closure issue. 
The matrix $\bK$ is usually large, as the state is lifted to a high-dimensional space, and the resulting model rarely closes and spurious eigenfunctions may appear.
Building on EDMD, we construct a reduced-order model using validated Koopman eigenfunctions, so that the resulting model will be \textit{closed by design}.}
Both KRONIC and EDMD with control~\cite{Korda2016arxiv} provide Koopman-based system identification that can be leveraged for model-based control, such as MPC or LQR, as shown in Fig.~\ref{Fig:KRONICvsEDMDc}.  
However, there are a number of key differences: 
(1) KRONIC directly identifies Koopman eigenfunctions, while EDMDc approximates the Koopman operator restricted to a high-dimensional span of observables.  
(2) EDMD augments the state vector with nonlinear measurements, increasing the dimension of the system.  In contrast, KRONIC yields a reduced-order model in terms of a few Koopman eigenfunctions.
(3) Control is incorporated in EDMDc as an approximated affine linear term. KRONIC derives an expression for how eigenfunctions are affected by control through the generator equation.  However, this may render the control term bilinear (nonlinear in the state). 
(4) The cost function for EDMDc is defined in the state or measurement space, while KRONIC defines the cost in eigenfunctions; note that these cost functions are not always transferable.  
(5) Finally, KRONIC readily admits more complicated solutions, such as limit cycle stabilization, as these correspond to level sets of the eigenfunctions.

As with previous studies, this work further cements the importance of accurate identification and representation of Koopman eigenfunctions.  
Future work will continue to develop algorithms to extract approximate eigenfunctions from data, and it is likely that these efforts will benefit from advances in machine learning. 
In addition, there is a fundamental uncertainty principle in representing Koopman eigenfunctions, as these eigenfunctions may themselves be irrepresentable, as are the long-time flow maps for chaotic systems. 
Instead of seeking perfect Koopman eigenfunctions, which may not be attainable, it will be important to incorporate uncertainty quantification into the data-driven Koopman framework.  
Model uncertainties may then be managed with robust control.  

The present work also highlights an important choice of perspective when working with Koopman approximations.  
Generally, Koopman eigenfunctions are global objects, such as the Hamiltonian energy function.  
Although a global, linear representation of the dynamics is appealing, there is also information that is stripped from these representations.  
For example, in the case of the Hamiltonian eigenfunction, information about specific fixed points and spatial locations are folded into a single scalar energy.  
If the Hamiltonian is viewed as a topography over the phase space, then this eigenfunction only carries information about the \emph{altitude}, and not the \emph{location}.  
In contrast, Lan and Mezi\'c~\cite{Lan2013physd} 
show that it is possible to extend the Hartman-Grobman theorem to the entire basin of attraction of certain fixed points and periodic orbits, providing a \emph{local} linear embedding of the dynamics.  
Connecting these perspectives will continue to yield interesting and important advances in Koopman theory.  
In addition, there are known connections between the eigenvalues of geometric structures in phase space and the spectrum of the Koopman operator~\cite{Mezic2017arxiv}.  
This knowledge may guide the accurate identification of Koopman eigenfunctions.

Formulating control in terms of Koopman eigenfunctions requires a change of perspective as the control objective may now be defined in eigenfunctions coordinates. Eigenfunctions characterize, e.g., geometric properties of the system, such as fixed points and limit cycles, as particular level sets, and the control objective can be equivalently formulated to steer the system towards these objects. 
Further, particular eigenfunctions represent coherent structures, i.e. persistent features of the system, which have been targeted in control for a long time. 
However, the specific selection of eigenfunctions to control and their interpretation regarding specific control goals remains an open problem.
Nevertheless, it may still be possible to formulate the control in the original state space, e.g. by incorporating the state as observables, modifying the state weight matrix appropriately, or by learning an approximation of the inverse mapping~\cite{kawahara2016nips}, which can be more easily incorporated in the context of model predictive control.
This also motivates additional work to understand how controllability and observability in these coordinates relate to properties of the nonlinear system. 
The degree of observability and controllability will generally vary with different eigenfunctions, so that it may be possible to obtain balanced realizations.  
Moreover, classic results, such as the PBH test, indicate that multi-channel actuation may be necessary to simultaneously control different eigenfunctions corresponding to the same eigenvalue, such as the Hamiltonian energy and conserved angular momentum.  
The additional degrees of freedom arising from multi-channel inputs can also be used for eigenstructure assignment to shape Koopman eigenfunctions~\cite{hemati2017aiaa}.  
Thus, actuation may modify both the shape of coherent structures (i.e., Koopman modes associated with a particular eigenfunction) and their time dynamics.  
It may also be possible to use Koopman linear embeddings to optimize sensor and actuator placement for nonlinear systems.  

Finally, as undamped or lightly damped eigenfunctions correspond to conserved or nearly conserved quantities, there are many potential applications of the proposed control strategy.  
For example, symmetries give rise to other conserved quantities, which will likewise yield new Koopman eigenfunctions. 
In many physical systems, simultaneously controlling the system energy and angular momentum may be an important goal.  
Much of the present work was formulated with the problem of space mission design in mind.  
Energy efficient transport throughout the solar system has long driven advances in theoretical and computational dynamical systems, and may stand to benefit from control based on Koopman eigenfunctions.  
More generally, there is a broad range of applications that stand to benefit from improved nonlinear control, include self-driving cars, the control of turbulence, suppressing the spread of disease, stabilizing financial markets, human machine interfaces, prosthetics and rehabilitation, and the treatment of neurological disorders, to name only a few.

%% file: ControlAppendix.tex
\RefereeOne{
We begin by formulating optimal control for the full nonlinear system in~\eqref{Eqn:OptimalControl:NonlinearSystemWithControl}.  
This procedure then simplifies considerably for linearly factorizable or linear dynamics, motivating Koopman embeddings.  
}

\RefereeOne{
\subsection{Hamilton-Jacobi-Bellman equation}
A nonlinear optimal control formulation for the system in~\eqref{Eqn:OptimalControl:NonlinearSystemWithControl} can be established using dynamic programming~\cite{Bertsekas2005book}. 
This relies on Bellman's principle of optimality~\cite{bellman1957book}, stating that an optimal trajectory remains optimal when the problem is initialized at intermediate steps, and leads to the Hamilton-Jacobi-Bellman (HJB) equation~\cite{bellman1964book}:
\begin{equation}\label{Eqn:OptimalControl:HJB}
0  = \inf\limits_{\bu}
\left\{ L(\bx,\bu)  + \nabla_{\bx} J^{*}\cdot{\boldsymbol{\it f}}(\bx,\bu)\right\} = \inf\limits_{\bu} H\left(\bx,\nabla_{\bx} J^{*},\bu\right).
\end{equation}
If the HJB equation~\eqref{Eqn:OptimalControl:HJB} has a continuously differentiable, positive definite solution $J^{*}(\bx)$~\cite{Bertsekas2005book}, also referred to as an optimal value function,
then the optimal control law is
\begin{equation}\label{Eqn:OptimalControl:HJB_ControlLaw}
\bu^{*} = \mathrm{arg}\,\min\limits_{\bu}\, H\left(\bx^{*},\frac{\partial J^{*}}{\partial x}, \bu \right). 
\end{equation}
The solution is a state-feedback control law, i.e.\ a \emph{closed-loop} controller, which is optimal for any initial condition and solved for all states at once.
Solving this nonlinear PDE is computationally challenging and only feasible for low-dimensional problems.
If solved, however, the HJB equation provides a global solution to the optimal control problem.
}

\RefereeOne{
\subsection{Euler-Lagrange equations}
An alternative approach searches for the optimal trajectory under certain constraints 
using a variational argument and Pontryagin's minimum (or maximum) principle~\cite{Pontryagin2062interscience}.
Here, a small variation around an optimal trajectory is considered, for which the change in cost shall vanish. 
The sequence of control inputs that solves this problem, constrained by the state dynamics~\eqref{Eqn:OptimalControl:NonlinearSystemWithControl}, constitutes the optimal control.
However, the Euler-Lagrange equations provide only necessary, but not sufficient, conditions, so that the solution may only converge to a local minimum and the control is generally not  optimal in practice.
}

\RefereeOne{
Similar to~\eqref{Eqn:OptimalControl:HJB}, a Lagrange-type expression is formulated for the control Hamiltonian $H(\bx,\bu,\bz) := L(\bx,\bu) + \bz^T {\boldsymbol{\it f}}(\bx,\bu)$, where $\bz$ is referred to as Lagrange multiplier, co-state or adjoint variable. 
A necessary condition is the minimization of the control Hamiltonian (i.e., Pontryagin's minimum principle), leading to the following set of coupled ordinary differential equations:
\begin{subequations}\label{Eq:OptimalControl:EulerLagrange}
\begin{align}
	\dot{\bx} =& \frac{\partial H^T}{\partial \bz}(\bx,\bu^{*},\bz) = {\boldsymbol{\it f}}(\bx,\bu)\\
	\dot{\bz} =&-\frac{\partial H^T}{\partial \bx}(\bx,\bu^{*},\bz)
\end{align}
\end{subequations}
with stationarity condition $\frac{\partial H^T}{\partial \bu}(\bx,\bu^{*},\bz)=0$ and transversality boundary conditions $\bx(0) = \bx_0$ and $\lim_{T\rightarrow\infty}\bz(T) = {\bf 0}$.
Then, the optimal control is given by
\begin{equation}
	\bu^{*} = \mathrm{arg}\,\min\limits_{\bu}\,H(\bx,\bu,\bz).
\end{equation}
If $L(\cdot)$ and ${\boldsymbol{\it f}}(\bx(\cdot),\bu(\cdot))$ are concave, then the necessary condition is also sufficient and any path that satisfies these conditions also solves the optimal control problem.
}

\RefereeOne{
The Euler-Lagrange equations~\eqref{Eq:OptimalControl:EulerLagrange} lead to an optimal \emph{open-loop} control: 
The control $\bu^{*}$ point-wise minimizes the control Hamiltonian $H(\bx^{*}(t),\bu,\bz^{*}(t))$.
This two-point boundary value (TPBV) problem is solvable for higher-dimensional problems in contrast to the HJB equation. 
However,
in practice, only a local minimum is found.
For systems with a high number of degrees of freedom, such as fluid flows, expensive direct and adjoint simulations render this approach infeasible, instead motivating the use of reduced-order models.
}

\RefereeOne{
\subsection{Linear or factorized systems}\label{Sec:Control:SDRE}
The control problem above simplifies considerably for linear systems of the form
\begin{equation}\label{Eqn:OptimalControl:LinearSystem}
\frac{\mathrm d}{\mathrm dt}\bx(t) = \bA \bx + \bB\bu,\quad \bx(0) = \bx_0.
\end{equation} 
Using \RefereeALL{ansatz $\bz = \bP\bx$} for the co-state, the optimal control is given by
\begin{equation}\label{Eqn:OptimalControl:Linear:ControlLaw}
\bu = -\bR^{-1}\bB^T\bz = -\bR^{-1}\bB^T\bP\bx,
\end{equation}
with constant gain $\bC = \bR^{-1}\bB^T\bP$ and
where the positive semi-definite matrix $\bP\in\mathbb{R}^{n\times n}$ is the solution to the algebraic Riccati equation (ARE)~\cite{stengel2012book}:
\begin{equation}\label{Eqn:OptimalControl:ARE}
\bQ + \bP\bA + \bA^T\bP - \bP\bB\bR^{-1}\bB^T\bP = {\bf 0}.
\end{equation}
For this special case, referred to as the linear quadratic regulator (LQR), the Euler-Lagrange equations and the HJB equation lead to the same solution, using a quadratic ansatz for co-state and value function.
The ARE~\eqref{Eqn:OptimalControl:ARE} can be solved upfront and yields a global state-feedback control law.
This simplicity motivates efforts to find linear representations for nonlinear systems and explains why nonlinear embeddings via Koopman operator theory are so appealing.
}

\RefereeOne{
For later reference, we point out an extension of LQR to control-affine, nonlinear systems by factoring the governing equations into a linear-like structure, as  in~\eqref{Eqn:OptimalControl:LinearSystem},
where the state transition matrix and actuation matrix become state-dependent:
\begin{equation}\label{Eqn:OptimalControl:FacorizedLinearSystem}
\frac{\mathrm d}{\mathrm dt}\bx(t) = \bA(\bx) \bx + \bB(\bx)\bu,\quad \bx(0) = \bx_0.
\end{equation}
For these systems, it is common to solve the ARE point-wise, i.e.\ at each point $\bx$ in time, leading to the state-dependent Ricatti (SDRE) equation~\cite{pearson1962ije} 
\begin{equation}\label{Eqn:OptimalControl:SDRE}
\bQ + \bP\bA(\bx) + \bA^T(\bx)\bP - \bP\bB(\bx)\bR^{-1}\bB^T(\bx)\bP = {\bf 0}
\end{equation}
where the weight matrices $\bQ$ and $\bR$ may also be state dependent.
In contrast to the ARE, which can be pre-computed off-line, the SDRE is solved on-line at each $\bx$. 
The required computational power and memory during the on-line phase may render it infeasible for high-dimensional systems.
The SDRE generalizes LQR for nonlinear systems, retaining a simple implementation and often yielding near optimal controllers~\cite{clautier1996proc}.  
We refer to~\cite{cimen2008ifac,cloutier1997ieee} for a review on the SDRE and related discussions, e.g.\ on controllability and stability.
}

%% file: ControlError.tex
\RefereeOne{
We examine the effect of an error in the representation of a Koopman eigenfunction on the closed-loop dynamics based on~\eqref{Eqn:KRONIC:SDRE} and provide an upper error bound. 
We consider the following cost function
\begin{equation}
	J = \int Q\varphi^2(\bx) + \bR \bu^2 \,\mathrm{dt}
\end{equation}	
where a single eigenfunction $\varphi(\bx)$ is to be controlled
and with positive weights $Q$ and $\bR$.
Assuming a control-affine system $\dot{\bx} = {\bf f}(\bx)+\bB(\bx)\bu$, 
the eigenfunction satisfies
\begin{equation}
	\dot{\varphi}(\bx) = \lambda\varphi(\bx) + \bC(\bx)u
\end{equation}	
with $\bC(\bx):=\nabla\varphi(\bx)\cdot \bB(\bx)$. Here, $\varphi$ is the  Koopman eigenfunction associated with the autonomous dynamics satisfying $\dot{\varphi} = \nabla\varphi(\bx)\cdot{\bf f}(\bx)=\lambda\varphi(\bx)$.
We assume full-state measurements $\bx$ are available and knowledge of how the control affects these measurements $\bB(\bx)$. Thus, we do not have to represent the control term $\bC(\bx)$ in terms of eigenfunctions but can evaluate it exactly given the current measurement $\bx$ and knowledge of $\varphi(\bx)$. 
Solving the scalar Riccati equation~\eqref{Eqn:KRONIC:SDRE} analytically yields
\begin{equation}
	P = \frac{1}{\bC(\bx)\bR^{-1}\bC^T(\bx)}\left(\lambda+\sqrt{\lambda^2+Q\bC(\bx)\bR^{-1}\bC^T(\bx)}\right).
\end{equation}	
The feedback control is then given by
\begin{equation}
	\bu = -\bR^{-1}\bC^TP\varphi = -\frac{\bR^{-1}\bC^T(\bx)}{\bC(\bx)\bR^{-1}\bC^T(\bx)} \left(\lambda+ \sqrt{\lambda^2+Q\bC(\bx)\bR^{-1}\bC^T(\bx)}\right) \varphi(\bx)
\end{equation}
and the resulting closed-loop dynamics are
\begin{equation}
	\dot{\varphi}(\bx) = -\sqrt{\lambda^2+Q\bC(\bx)\bR^{-1}\bC^T(\bx)}\varphi(\bx) = -\sqrt{\mu}\varphi(\bx),
\end{equation}	
where $\mu := \lambda^2+Q\bC(\bx)\bR^{-1}\bC^T(\bx)$. We examine the effect of a misrepresentation of the eigenfunction:
\begin{equation}
	\hat{\varphi}(\bx):= \varphi(\bx) + \varepsilon \psi(\bx)
\end{equation}
where $\varepsilon \psi(\bx)$ is the discrepancy to the true eigenfunction $\varphi(\bx)$ with small $\varepsilon$.
The dynamics of $\hat{\varphi}(\bx)$ are given in a factorized representation by
\begin{eqnarray}
	\dot{\hat{\varphi}}(\bx) = \dot{\varphi}(\bx) + \varepsilon\dot{\psi}(\bx)
	&=& \lambda\varphi(\bx) + \bC(\bx)\bu + \varepsilon \nabla\psi\cdot({\bf f}(\bx)+\bB(\bx)\bu)\notag\\
	&=& \frac{\lambda\varphi(\bx)+\varepsilon\nabla\psi(\bx)\cdot{\bf f}(\bx)}{\varphi(\bx)+\varepsilon\psi(\bx)}\hat{\varphi}(\bx) + \left(\bC(\bx) + \varepsilon\bD(\bx) \right)\bu
\end{eqnarray}
with $\bD(\bx):=\nabla\psi(\bx)\cdot\bB(\bx)$. 
The upper bound for the error in $\mu$ due to the misrepresentation of $\varphi(\bx)$ is
\begin{eqnarray}
	\vert\mu-\hat{\mu} \vert 
	&=& \Big\lvert\lambda^2+Q\bC(\bx)\bR^{-1}\bC^T(\bx) - \left(  \frac{\lambda\varphi(\bx)+\varepsilon\nabla\psi(\bx)\cdot{\bf f}(\bx)}{\varphi(\bx)+\varepsilon\psi(\bx)} \right)^2
	-Q(\bC(\bx)+\varepsilon\bD(\bx))\bR^{-1}(\bC(\bx)+\varepsilon\bD(\bx))^T\Big\rvert\notag\\
	&=& \Big\lvert\lambda^2-\varepsilon Q (\bC\bR^{-1}\bD^T+\bD\bR^{-1}\bC^T)-\varepsilon^2Q\bD\bR^{-1}\bD^T-\frac{\lambda^2\varphi^2 + \varepsilon 2\lambda\varphi\nabla\psi\cdot{\bf f}+\varepsilon^2(\nabla\psi\cdot{\bf f})^2}{(\varphi+\varepsilon\psi)^2}\Big\rvert\notag\\
	&\leq& \Big\lvert\lambda^2-\varepsilon Q (\bC\bR^{-1}\bD^T+\bD\bR^{-1}\bC^T)-\varepsilon^2Q\bD\bR^{-1}\bD^T-\frac{\lambda^2\varphi^2 + \varepsilon 2\lambda\varphi\nabla\psi\cdot{\bf f}+\varepsilon^2(\nabla\psi\cdot{\bf f})^2}{\varphi^2(\bx)}\Big\rvert\notag\\
	&=& \Big\lvert -\varepsilon Q (\bC\bR^{-1}\bD^T+\bD\bR^{-1}\bC^T)-\varepsilon^2Q\bD\bR^{-1}\bD^T-\frac{\varepsilon 2\lambda\varphi\nabla\psi\cdot{\bf f}+\varepsilon^2(\nabla\psi\cdot{\bf f})^2}{\varphi^2}\Big\rvert\notag\\
	&=& \Big\lvert-\varepsilon \left(Q\bC\bR^{-1}\bD^T+Q\bD\bR^{-1}\bC^T+2\lambda\frac{\nabla\psi\cdot{\bf f}}{\varphi(\bx)}\right)-\varepsilon^2\left(Q\bD\bR^{-1}\bD^T+\frac{(\nabla\psi\cdot{\bf f})^2}{\varphi^2(\bx)}\right)\Big\rvert.
\end{eqnarray}
}	

%% file: AnalyticalEigenfunction.tex
It is possible to solve the PDE in \eqref{Eq:KoopmanPDE} for eigenfunctions using standard techniques, such as Taylor or Laurent series.  

\subsection{Linear dynamics}
Consider linear dynamics
\begin{equation}
\frac{d}{dt}x =  x.
\end{equation}
Assuming a Taylor series expansion for $\varphi(x)$:
\begin{align*}
\varphi(x) &= c_0 + c_1x + c_2x^2 + c_3x^3 + \cdots
\end{align*}
then the gradient and directional derivatives are given by:
\begin{align*}
\nabla \varphi & = c_1 + 2c_2x + 3 c_3x^2 + 4c_4x^3+\cdots\\
\nabla\varphi\cdot f & = c_1 x + 2c_2 x^2 +3c_3 x^3 + 4c_4 x^4+\cdots
\end{align*}
Solving for terms in \eqref{Eq:KoopmanPDE}, $c_1=0$ must hold.  
For any positive integer $\lambda$, only one coefficient may be nonzero.  
Specifically, for $\lambda=k\in\mathbb{Z}^+$, then $\varphi(x) = cx^k$ is an eigenfunction for any constant $c$. 
For instance, if $\lambda=1$, then $\varphi(x) =  x$.  

\subsection{Quadratic nonlinear dynamics}
Consider a nonlinear dynamical system
\begin{align}
\frac{d}{dt}=x^2.
\end{align}
There is no Taylor series that satisfies \eqref{Eq:KoopmanPDE}, except the trivial solution $\varphi=0$ for $\lambda=0$.  
Instead, we use a Laurent series:
\begin{align*}
\varphi(x) &= \cdots +c_{-3}x^{-3} + c_{-2}x^{-2} + c_{-1}x^{-1} +  c_0 \notag\\
& \hspace{0.9cm}+c_1x + c_2x^2 + c_3x^3 + \cdots.
\end{align*}
The gradient and directional derivatives are given by:
\begin{align*}
\nabla \varphi & =\cdots -3c_{-3}x^{-4} - 2c_{-2}x^{-3} -c_{-1}x^{-2} + c_1 + 2c_2x \notag\\
					  &\hspace{0.4cm}+ 3 c_3x^2 + 4c_4x^3+\cdots\\
\nabla\varphi\cdot f & =\cdots -3c_{-3}x^{-2} - 2c_{-2}x^{-1} -c_{-1} + c_1x^2 + 2c_2x^3 \notag\\
					 &\hspace{0.4cm}+ 3 c_3x^4 + 4c_4x^5+\cdots.
\end{align*}
All coefficients with positive index are zero, i.e. $c_k=0$ for $k\geq 1$.  
The nonpositive index coefficients are given by
 the recursion $\lambda c_{k+1} = k c_{k}$, for $k\leq -1$.  
%
The Laurent series is
\begin{align*}
\varphi(x) &= c_0\left(1 - \lambda x^{-1} + \frac{\lambda^2}{2}x^{-2} - \frac{\lambda^3}{3!}x^{-3}+\cdots\right)=c_0e^{-\lambda/x}.
\end{align*}
This holds for all $\lambda\in\mathbb{C}$.  There are also other Koopman eigenfunctions that can be identified from the series.

\subsection{Polynomial nonlinear dynamics}
For a more general nonlinear dynamical system
\begin{align}
\frac{d}{dt}=ax^n,
\end{align}
$\varphi(x)=e^{\frac{\lambda}{(1-n)a}x^{1-n}}$ is an eigenfunction for all $\lambda\in\mathbb{C}$.  